\documentclass[a4paper,12pt]{amsart}
\usepackage{amsfonts}
\usepackage{amssymb}
\usepackage{ifthen}
\usepackage{amscd}
\usepackage{amsxtra}
\usepackage{graphicx}
\usepackage{color}
\nonstopmode \numberwithin{equation}{section}
\newtheorem{thm}{Theorem}%[section]
\newtheorem{lem}{Lemma}%[section]
\newtheorem{cor}{Corollary}%[section]
\newtheorem{cl}{Claim}%[section]
\newtheorem{ca}{Case}%[section]
\newtheorem{sca}{Subcase}%[section]
\newtheorem{scl}{Subclaim}%[section]
\newtheorem{conj}{Conjecture}

\theoremstyle{definition}
\newtheorem{defn}{Definition}%[section]

\newtheorem{op}[equation]{Open Problem}
\newtheorem{ques}[equation]{Question}
\newtheorem{rem}{Remark}%[section]
\newtheorem{exam}[equation]{Example}

\newcounter {own}
\def\theown {\thesection       .\arabic{own}}

\newenvironment{pf}[1][]{%
 \vskip 3mm
 \noindent
 \ifthenelse{\equal{#1}{}}%
  {{\slshape Proof. }}%
  {{\slshape #1.} }%
 }%
{\qed\bigskip}

\newcounter{alphabet}
\newcounter{tmp}
\newenvironment{Thm}[1][]{\refstepcounter{alphabet}%
\bigskip%
\noindent%
{\bf Theorem \Alph{alphabet}}%
\ifthenelse{\equal{#1}{}}{}{ (#1)}%
{\bf .} \itshape}{\vskip 8pt}

\makeatletter
\newcommand{\Ref}[1]{\@ifundefined{r@#1}{}{\setcounter{tmp}{\ref{#1}}\Alph{tmp}}}
\makeatother
\newenvironment{Lem}[1][]{\refstepcounter{alphabet}%
\bigskip%
\noindent%
{\bf Lemma \Alph{alphabet}}%
{\bf .} \itshape}{\vskip 8pt}

\newcommand{\IR}{{\mathbb R}}

\newcommand{\ID}{{\mathbb D}}

\newcommand{\diam}{{\operatorname{diam}}}

\def\be{\begin{equation}}
\def\ee{\end{equation}}

\newcommand{\bee}{\begin{enumerate}}
\newcommand{\eee}{\end{enumerate}}

\newcommand{\blem}{\begin{lem}}
\newcommand{\elem}{\end{lem}}
\newcommand{\bthm}{\begin{thm}}
\newcommand{\ethm}{\end{thm}}
\newcommand{\bcor}{\begin{cor}}
\newcommand{\ecor}{\end{cor}}
\newcommand{\beg}{\begin{exam}}
\newcommand{\eeg}{\end{exam}}
\newcommand{\begs}{\begin{examples}}
\newcommand{\eegs}{\end{examples}}
\newcommand{\bdefe}{\begin{defn}}
\newcommand{\edefe}{\end{defn}}
\newcommand{\bprob}{\begin{prob}}
\newcommand{\eprob}{\end{prob}}
\newcommand{\bques}{\begin{ques}}
\newcommand{\eques}{\end{ques}}
\newcommand{\bei}{\begin{itemize}}
\newcommand{\eei}{\end{itemize}}
\newcommand{\bcon}{\begin{conj}}
\newcommand{\econ}{\end{conj}}
\newcommand{\bop}{\begin{op}}
\newcommand{\eop}{\end{op}}

\newcommand{\bca}{\begin{ca}}
\newcommand{\eca}{\end{ca}}
\newcommand{\bsca}{\begin{sca}}
\newcommand{\esca}{\end{sca}}

\newcommand{\bcl}{\begin{cl}}
\newcommand{\ecl}{\end{cl}}

\newcommand{\bscl}{\begin{scl}}
\newcommand{\escl}{\end{scl}}

\newcommand{\bcons}{\begin{conjs}}
\newcommand{\econs}{\end{conjs}}
\newcommand{\bprop}{\begin{propo}}
\newcommand{\eprop}{\end{propo}}
\newcommand{\br}{\begin{rem}}
\newcommand{\er}{\end{rem}}
\newcommand{\brs}{\begin{rems}}
\newcommand{\ers}{\end{rems}}
\newcommand{\bo}{\begin{obser}}
\newcommand{\eo}{\end{obser}}
\newcommand{\bos}{\begin{obsers}}
\newcommand{\eos}{\end{obsers}}
\newcommand{\bpf}{\begin{pf}}
\newcommand{\epf}{\end{pf}}
\newcommand{\ba}{\begin{array}}
\newcommand{\ea}{\end{array}}
\newcommand{\beq}{\begin{eqnarray}}
\newcommand{\beqq}{\begin{eqnarray*}}
\newcommand{\eeq}{\end{eqnarray}}
\newcommand{\eeqq}{\end{eqnarray*}}

\newcommand{\ds}{\displaystyle}

\newcounter{minutes}
\divide\time by 60
\newcounter{hours}
\multiply\time by 60 \addtocounter{minutes}{-\time}
\begin{document}

\bibliographystyle{amsplain}
\title []
{Radial length, radial John disks and $K$-quasiconformal harmonic
mappings }

%%%%%%%% BEGIN TIMESTAMP
\def\thefootnote{}
\footnotetext{ \texttt{\tiny File:~\jobname .tex,
          printed: \number\day-\number\month-\number\year,
          \thehours.\ifnum\theminutes<10{0}\fi\theminutes}
} \makeatletter\def\thefootnote{\@arabic\c@footnote}\makeatother
%%%%%%%% END TIMESTAMP

\author{Shaolin Chen}
\address{S. L. Chen, College of Mathematics and Statistics, Hengyang Normal University, Hengyang, Hunan 421008,
People's Republic of China.} \email{mathechen@126.com}

\author{Saminathan Ponnusamy %$^\dagger $ %${}^{~\mathbf{*}}$
}
\address{S. Ponnusamy,
Indian Statistical Institute (ISI), Chennai Centre, SETS (Society
for Electronic Transactions and Security), MGR Knowledge City, CIT
Campus, Taramani, Chennai 600 113, India. }
\email{samy@isichennai.res.in, samy@iitm.ac.in}

%\author{ A. Rasila }
%\address{A. Rasila, Department of Mathematics and Systems Analysis, Aalto University, P. O. Box 11100, FI-00076 Aalto,
% Finland.} \email{antti.rasila@iki.fi}

%\author{X. Wang$^{\mathbf{*}}$
%${}^{~\mathbf{*}}$
%}
%\address{X. Wang, Department of Mathematics,
%Hunan Normal University, Changsha, Hunan 410081, People's Republic
%of China.} \email{xtwang@hunnu.edu.cn}

 \keywords{$K$-quasiconformal harmonic mapping, radial John
disk, radial length,
 Pommerenke interior domain.
\\
%$%{}^{\mathbf{*}}
%^\dagger${\tt ~~%Corresponding author.
% This author is on leave from
%the Department of Mathematics, Indian Institute of Technology
%Madras, Chennai-600 036, India}
%%\\ ${}^{\mathbf{*}}$ Corresponding author
This second author is on leave from IIT Madras.
}
%\date{\today } %November 4, 10;
%File: Ch-W-S12${}_{}$equiv-mod${}_{}$submit.tex}

\begin{abstract}
In this article, we continue our investigations of
the boundary behavior of harmonic mappings. We first discuss the
classical problem on the growth of radial length and obtain a sharp
growth theorem of the radial length of $K$-quasiconformal harmonic
mappings. Then we present an alternate characterization of radial John
disks. In addition, we investigate the linear measure distortion and the
Lipschitz continuity on $K$-quasiconformal harmonic mappings of the
unit disk onto a radial John disk. Finally, using  Pommerenke interior domains,
we characterize certain differential properties of $K$-quasiconformal harmonic mappings

\end{abstract}

%\thanks{The research was partly supported by
%NSF of China (No. 11071063)} %and  Hunan Provincial Innovation
%%Foundation for Postgraduate (No. 125000-4113).  }

\maketitle \pagestyle{myheadings} \markboth{ S. L. Chen and S.
Ponnusamy}{Radial length, radial John disks and $K$-quasiconformal
harmonic mappings}

\section{Introduction and statement of main results }\label{csw-sec1}

%For $a\in\mathbb{C}$ and  $r>0$, we let $\ID(a,r)=\{z:\, |z-a|<r\}$
%so that $\mathbb{D}_r:=\mathbb{D}(0,r)$ and thus,
%$\mathbb{D}:=\mathbb{D}_1$ denotes the open unit disk in the complex
%plane $\mathbb{C}$.

This paper continues the study of previous work of the authors
\cite{CP} and is mainly motivated by the articles of Beardon and Carne \cite{BC},
Carroll and Twomey \cite{CT},  Chuaqui et al.  \cite{COP}, Pommerenke \cite{Po},
and the monograph of Pommerenke \cite{Po1}.
%The present article is organized as follows. We first investigate the growth
%of the {\it radial length} of $K$-quasiconformal harmonic mappings
%(see Theorem \ref{thm-CPx}). Then we give a new
%characterization of the radial John disk (see Theorem \ref{thm-J}).
%and the Lipschitz continuity on  $K$-quasiconformal harmonic
%mappings of the unit disk $|z|<1$ onto {\it radial John disks} (see Theorems
%\ref{thm-c1} and \ref{thm-c3}). Some differential properties of the
%$K$-quasiconformal harmonic mappings are also characterized by
%using {\it Pommerenke interior domains} (see Theorems \ref{thm-1},
%\ref{thm-1.02} and Corollary \ref{cor-1}).
In order to state our first result concerning the growth of the {\it radial length} of $K$-quasiconformal
harmonic mappings (see Theorem \ref{thm-CPx}), we need to recall some basic definitions and some results which
motivate the present work.

Let $f$ be a complex-valued and continuously differentiable function
defined in the unit disk $\ID=\{z:\, |z|<1\}$ and let
$\ell_{f}(\theta,r)$ be the length of the $f$-image (with counting
multiplicity) of the  {\it radial line segment} $[0,z]$ from $0$ to
$z=re^{i\theta}\in \ID$, where $\theta\in[0,2\pi]$ is fixed and
$r\in[0,1)$. Then (cf. \cite{CLP})
$$\ell_{f}(\theta,r):=\ell\big(f([0, z])\big)=\int_{0}^{r}\left |\,df(\rho
e^{i\theta})\right |=\int_{0}^{r}\left |f_{z}(\rho
e^{i\theta})+e^{-2i\theta}f_{\overline{z}}(\rho e^{i\theta})\right
|\,d\rho.
$$
In \cite{Ke}, Keogh showed that if $f$ is a bounded, analytic and
univalent function in $\mathbb{D}$, then, for each $\theta\in[0,2\pi]$,
\be\label{eq-cpk1} \ell_{f}(\theta,r)=O \left (
\big(\log (1/(1-r))\big)^{1/2}\right) ~\mbox{ as $r\rightarrow
1^{-}$}.
\ee
Throughout the discussion, we let
\be\label{CP-extraeq2} \psi(r)=\big(\log (1/(1-r) )\big)^{1/2}~\mbox{ for $0<r<1$}.
\ee
Keogh also gave some examples to show that
the exponent $1/2$ in (\ref{eq-cpk1}) can not be decreased. Jenkins
improved on these examples in \cite{Je}, and Kennedy \cite{Ken}
presented further examples by showing that
%completed this sequence of results for bounded and univalent analytic functions
%We begin  by recalling some known results relating to the growth of
%$\ell_{f}(\theta,r)$ for bounded and univalent analytic functions.
%\be\label{eq-cpk2}
$$\ell_{f}(\theta,r)=O (\mu(r)\psi(r) )~\mbox{as $r\rightarrow 1^{-}$}
$$
is false in general for every positive function
$\mu$ in $[0,1)$ satisfying $\mu(r)\rightarrow0$ as
$r\rightarrow1^{-}$. In \cite{CT}, Carroll and Twomey established certain refinements and extension of these results
without the boundedness condition in the following form.

\begin{Thm}\label{Thm-CT}
Suppose that $f(z)=a_{1}z+a_{2}z^{2}+\cdots$ is univalent in
$\mathbb{D}$. Then, for  any fixed $\theta\in[0,2\pi]$, there is a
constant $C_{1}>0$ such that
\be\label{CP-extraeq1}
\ell_{f}(\theta,r)\leq C_{1}\max_{\rho\in[0,r]}|f(\rho e^{i\theta})| \psi(r) %\left(\log\frac{1}{1-r}\right)^{1/2}
~\mbox{ for $r\in (0.5,1)$}. \ee If, further, $f(re^{i\theta})=O(1)$
as $r\rightarrow1^{-},$ then {\rm (\ref{eq-cpk1})} holds.
\end{Thm}

Later, Beardon and Carne \cite{BC} gave a relatively simple argument to Theorem \Ref{Thm-CT} in hyperbolic
geometry and provided with further examples. It is worth pointing out here two results
which strengthened \eqref{CP-extraeq1} and was inspired by the work of Sheil-Small \cite{Small69}
and Hall \cite{Hall76}. If $f\in {\mathcal S}$ is starlike, i.e. $f(\ID)$ contains the
line segment $[0,w]$ whenever it contains $w$, then (see \cite{Karu83})
$$\ell_{f}(\theta,r)\leq  |f(re^{i\theta})|  (1+r)<2 |f(re^{i\theta})|~\mbox{ for $r\in (0,1)$}
$$
and the inequality of course is not sharp for all $r$, but the bound $2$ sharp as the Koebe function
$k(z)=z/(1-z)^2$ shows and is attained when $r$ approaches $1$ (see \cite{Hall76,Small69}).
Later in 1993,  Balasubramanian et al. \cite{BKP93} showed that if $f\in {\mathcal S}$ is convex,
i.e. $f(\ID)$ is a convex domain, then
$$\ell_{f}(\theta,r)\leq  |f(re^{i\theta})|\, r^{-1} \arcsin r ~\mbox{ for $r\in (0,1)$}
$$
and the inequality is sharp as the convex function $f(z)=z/(1-z)$
shows. Note that $\varphi (r)=r^{-1} \arcsin r$ is increasing on
$(0,1)$ and $\varphi (r)\leq \lim_{r\rightarrow 1^{-}}\varphi (r)
=\pi /2$ and thus, the conjecture of Hall \cite{Hall80} was settled
(see also \cite{BP96}).

The first  aims of this paper is to extend Theorem \Ref{Thm-CT} for the
case of harmonic quasiconformal mappings (see Theorem \ref{thm-CPx} below).
We need some preparation to state this result.

For a real $2\times2$ matrix $A$, we use the matrix norm
$\|A\|=\sup\{|Az|:\,|z|=1\}$ and the matrix function
$l(A)=\inf\{|Az|:\,|z|=1\}$.
For $z=x+iy\in\mathbb{C}$, the formal derivative of the complex-valued function $f=u+iv$ is given by the
Jacobian matrix
$$D_{f}=\left(\begin{array}{cccc}
\ds u_{x}\;~~ u_{y}\\[2mm]
\ds v_{x}\;~~ v_{y}
\end{array}\right),
$$
so that
$$\|D_{f}\|=|f_{z}|+|f_{\overline{z}}| ~\mbox{ and }~ l(D_{f})=\big| |f_{z}|-|f_{\overline{z}}|\big |,
$$
where $f_{z}=(1/2)\big( f_x-if_y\big)$ and $f_{\overline{z}}=(1/2)\big(f_x+if_y\big)$. Let $\Omega$ be a domain
in $\mathbb{C}$, with non-empty boundary. A sense-preserving
homeomorphism $f$ from a domain $\Omega$ onto $\Omega'$, contained
in the Sobolev class $W_{loc}^{1,2}(\Omega)$, is said to be a
{\it $K$-quasiconformal mapping} if, for $z\in\Omega$,
$$\|D_{f}(z)\|^{2}\leq K\big | \det D_{f}(z)\big |,~\mbox{i.e.,}~\|D_{f}(z)\|\leq Kl\big(D_{f}(z)\big),
$$
where $K\geq1$ and $\det D_{f}$ is the determinant of $D_{f}$ (cf. \cite{K,LV,Va,V}).

Let ${\mathcal S}_{H}$ denote the family of sense-preserving planar harmonic univalent mappings
$f=h+\overline{g}$ in $\mathbb{D}$, with the normalization $h(0)=g(0)=0$ and $h'(0)=1$.
Recall that $f$ is sense-preserving if the Jacobian $J_{f}$ of $f$ given by
$$J_{f}:=\det D_{f} =|f_{z}|^{2}-|f_{\overline{z}}|^{2} =|h'|^2-|g'|^2
$$
is positive. Thus, $f$ is locally univalent and sense-preserving in $\ID$ if and only if
$J_{f}(z)>0$ in $\ID$; or equivalently if $h'\neq 0$ in $\ID$ and the dilatation
$\omega =g'/h'$ has the property that $|\omega (z)|<1$ in $\ID$ (see \cite{Clunie-Small-84,Du,Lewy}).
The  family ${\mathcal S}_{H}$ together with a few other geometric
subclasses, originally investigated in detail by \cite{Clunie-Small-84,Small}, became instrumental
in the study of univalent harmonic mappings (see \cite{Du,SaRa2013}) and has attracted the attention of many
function theorists. If the co-analytic part
$g$ is identically zero in the decomposition of $f=h+\overline{g}$, then the class
${\mathcal S}_{H}$ reduces to the classical family $\mathcal S$ of
all normalized analytic univalent functions
$h(z)=z+\sum_{n=2}^{\infty}a_{n}z^{n}$ in $\ID$. If ${\mathcal
S}_H^{0}=\{f=h+\overline{g} \in {\mathcal S}_H: \,g'(0)=0 \} $, then
the family ${\mathcal S}_H^{0}$ is both normal and compact. See \cite{Clunie-Small-84} and also
\cite{CPRW,CP,Du,SaRa2013}.
%Recall that $f$ is sense-preserving if the Jacobian $J_{f}$ of $f$ is given by
%$$J_{f}:=\det D_{f} =|f_{z}|^{2}-|f_{\overline{z}}|^{2} =|h'|^2-|g'|^2
%$$
%is positive. Thus, $f$ is locally univalent and sense-preserving in $\ID$ if and only if
%$J_{f}>0$ in $\ID$; or equivalently if $h'\neq 0$ in $\ID$ and the dilatation
%$\omega =g'/h'$ has the property that $|\omega|<1$ in $\ID$ (see \cite{Clunie-Small-84,Du,Lewy}).

\begin{thm}\label{thm-CPx}
For $K\geq1$, let $f\in {\mathcal S}_{H}$ be a $K$-quasiconformal
harmonic mapping. Then, for any fixed $\theta\in[0,2\pi]$, there is
a constant $C_{2}>0$ such that
$$\ell_{f}(\theta,r)\leq C_{2}\max_{\rho\in[0,r]}|f(\rho e^{i\theta})| \psi(r) ~\mbox{ for $r\in (0.5,1)$}.
%\left(\log\frac{1}{1-r}\right)^{1/2}~\mbox{for}~r\in(\frac{1}{2},1).
$$
If, further, $f(re^{i\theta})=O(1)$ as $r\rightarrow1^{-},$ then
$$\ell_{f}(\theta,r)= O (\psi(r))~\mbox{ as $r\rightarrow1^{-}$},
$$
and the exponent $1/2$ in $\psi(r)$ defined by \eqref{CP-extraeq2}
cannot be replaced by a smaller number.
% is best possible.
\end{thm}

First we remark that if $K=1$, then Theorem \ref{thm-CPx} coincides with Theorem \Ref{Thm-CT}. Secondly,
the proof of Theorem \ref{thm-CPx} is substantially harder than the proof
of Theorem \Ref{Thm-CT}. This is because Beardon and Carne's argument of Theorem \Ref{Thm-CT} in
\cite{BC} is not applicable in the proof of Theorem \ref{thm-CPx}.

We need further notation and terminology before stating our second result.
Let $d_{\Omega}(z)$ be the Euclidean distance from $z$ to the boundary $\partial \Omega$ of $\Omega$. If $\Omega =\ID$, then we set
$d(z):=d_{\ID}(z)$.

%In particular, we always use
%$d(z)$ to denote the Euclidean distance from $z$ to the boundary
%$\partial \mathbb{D}$ of $\mathbb{D}.$

\begin{defn}\label{CP-de-1}
A bounded simply connected plane domain $G$ is called a {\it
$c$-John disk} for $c\geq1$ with {\it John center} $w_{0}\in G$ if
for each $w_{1}\in G$ there is a rectifiable arc $\gamma$, called a
{\it John curve}, in $G$ with end points $w_{1}$ and $w_{0}$ such that
%\be\label{eq-1}
$$\sigma_{\ell}(w)\leq cd_{G}(w)
$$
for all $w$ on $\gamma$, where $\gamma[w_{1},w]$ is the subarc of
$\gamma$ between $w_{1}$ and $w$, and $\sigma_{\ell}(w)$ is the
Euclidean length of $\gamma[w_{1},w]$ (see \cite{CP,KH,John,NV,Po1}).
\end{defn}

\begin{rem}
If $f$ is a complex-valued and univalent mapping in $\mathbb{D}$,
$G=f(\mathbb{D})$ and, for $z\in\mathbb{D}$, $\gamma=f([0,z])$ in
Definition \ref{CP-de-1}, then we call $c$-John disk  a {\it radial}
$c$-John disk, where $w_{0}=f(0)$ and $w=f(z)$. In particular, if
$f$ is a conformal mapping, then we call $c$-John disk  a {\it
hyperbolic} $c$-John disk.  It is well known that any point
$w_{0}\in G$ can be chosen as a John center by modifying the
constant $c$ if necessary. When we do not wish to emphasize the role of $c$,
then we regard the $c$-John disk simply as a John disk in the natural way
(cf. \cite{CP,KH,John,NV}).
\end{rem}

Unless otherwise stated, throughout the discussion we consider the
following terminology. Denote by ${\mathcal F}(K)$ if $f\in
{\mathcal F}$ and is a $K$-quasiconformal harmonic mapping in $\ID$,
where $K\geq 1$. Also, we denote by ${\mathcal F}(K,\Omega)$ if
$f\in {\mathcal F}(K)$ and $f$ maps $\ID$ onto $\Omega$. We prove
several results mainly when ${\mathcal F}$ equals one of ${\mathcal
S}_{H}$, ${\mathcal S}_{H}^{0}$, and ${\mathcal S}_{H_{2}}$, and
$\Omega$ equals either radial John disk or Pommerenke interior
domain.

Further, for $z\in\mathbb{D},$ we define
\be\label{CP-extraeq3}
B(z):=\{\zeta:\,|z|\leq|\zeta|<1,~|\arg z-\arg \zeta|\leq\pi(1-|z|)\}.
\ee
In the following, we continue our previous work of \cite{CP}, and give an another characterization of the radial John disk.

\begin{thm}\label{thm-J}
%For $K\geq1$, let $f\in{\mathcal S}_{H}^{0}$ be a $K$-quasiconformal harmonic mapping.

Let $f\in{\mathcal S}_{H}^{0}(K)$. Then the following are equivalent:
\begin{enumerate}
\item[{\rm (i)}]  $\Omega=f(\mathbb{D})$ is a radial John disk.

\item[{\rm (ii)}] There is an $x\in(0,1)$ such that
$$\sup_{|\zeta|=1}\sup_{r\in(0,1)}\frac{(1-\rho^{2})\|D_{f}(\rho\zeta)\|}
{(1-r^{2})\|D_{f}(r\zeta)\|}<1~\mbox{for}~\rho=\frac{x+r}{1+xr}.
$$

\item[{\rm (iii)}] $\ds \sup_{z\in\mathbb{D},~w\in B(z)}\frac{|f(z)-f(w)|}{(1-|z|^{2})\|D_{f}(z)\|}<\infty.$
\end{enumerate}
%Then $\Omega=f(\mathbb{D})$ is a radial John disk if and only if
%there is an $x\in(0,1)$ such that
%$$\sup_{|\zeta|=1}\sup_{r\in(0,1)}\frac{(1-\rho^{2})\|D_{f}(\rho\zeta)\|}
%{(1-r^{2})\|D_{f}(r\zeta)\|}<1~\mbox{for}~\rho=\frac{x+r}{1+xr}.$$
\end{thm}

Next, we  establish the linear measure distortion on $K$-quasiconformal harmonic mappings of $\mathbb{D}$
onto a radial John disk.

\begin{thm}\label{thm-c1}
%For $K\geq1$, let $f=h+\overline{g}\in{\mathcal S}_H^{0}$ be a
%$K$-quasiconformal harmonic mapping from $\mathbb{D}$ onto a radial
%John disk $\Omega$, where $h$ and $g$ are analytic in $\mathbb{D}$.
Let $f=h+\overline{g}\in{\mathcal S}_H^{0}(K,\Omega)$, where $\Omega$ is a radial John disk.
Then, for $a_{1},a_{2}\in\mathbb{D}$ with $B(a_{1})\subset
B(a_{2})$, there is a positive constant $C_{3}$ such that
$$\frac{\diam f(B(a_{1}))}{\diam f(B(a_{2}))}\leq C_{3}\left(\frac{\ell(B(a_{1})\cap\partial\mathbb{D})}
{\ell(B(a_{2})\cap\partial\mathbb{D})}\right)^{\alpha}.
$$
where $\alpha=\sup_{f\in{\mathcal S}_H} \frac{|h''(0)|}{2}$ and $B(z)$ is defined by \eqref{CP-extraeq3}.
\end{thm}

We remark that $2\leq\alpha=\sup_{f\in{\mathcal S}_{H}}\frac{|h''(0)|}{2}<\infty$, but the sharp value of $\alpha$
is still unknown (see \cite{CP,CHM,Du,Small}). We discuss the
Lipschitz continuity on $K$-quasiconformal harmonic mappings of
$\mathbb{D}$ onto a radial John disk, which is as follows.

\begin{thm}\label{thm-c3}
%For $K\geq1$, let $f\in{\mathcal S}_H^{0}$ be a $K$-quasiconformal
%harmonic mapping from $\mathbb{D}$ onto a radial John disk $\Omega$.
Let $f=h+\overline{g}\in{\mathcal S}_H^{0}(K,\Omega)$, where $\Omega$ is a radial John disk.
Then, for $z\in\mathbb{D}$ with $|z|\geq\frac{1}{2}$ and $\zeta_{1},
\zeta_{2}\in B(z)$,  there are constants $\delta_{1}\in(0,1)$ and
$C_{4}>0$ such that
$$|f(\zeta_{1})-f(\zeta_{2})|\leq C_{4}d_{\Omega}(f(z))\left(\frac{|\zeta_{1}-\zeta_{2}|}{1-|z|}\right)^{\delta_{1}}.
$$
\end{thm}

%Let $f\in{\mathcal S}_{H}$ be a $K$-quasiconformal harmonic mapping from $\mathbb{D}$ onto a domain $G$.

Let $f\in{\mathcal S}_H(K,G)$, where $G$ is domain. For $0<r<1$, let $\mathbb{D}_{r}=\{z:\, |z|<r\}$ and $\partial\mathbb{D}_{r}$ denote the
boundary of $\mathbb{D}_{r}$. Now, for  $w_{1}, w_{2}\in f\left(\partial\mathbb{D}_{r}\right)$, let $\gamma_{r}$ be
the smaller subarc of $f\left(\partial\mathbb{D}_{r}\right)$ between $w_{1}$ and $w_{2}$, and let
$$d_{G_{r}}(w_{1},w_{2})=\inf_{\Gamma}\diam \Gamma,
$$
where $\Gamma$ runs through all arcs from $w_{1}$ to $w_{2}$ that
lie in $G_{r}=f(\mathbb{D}_{r})$ except for their endpoints. If
\be\label{eq-y}
\sup_{0<r<1}\left\{\sup_{w_{1},w_{2}\in\gamma_{r}}\frac{\ell\big(\gamma_{r}[w_{1},w_{2}]\big)}{d_{G_{r}}(w_{1},w_{2})}\right\}<\infty,
\ee
then we call $G$  a {\it Pommerenke interior domain} (cf. \cite{CP,Po}). In particular, if $G$ is bounded, then we call $G$ as
a {\it bounded Pommerenke interior domain}.

Given a sense-preserving harmonic mapping $f=h+\overline{g}$ in $\ID$, fix $\zeta\in\ID$ and perform a disk automorphism
(also called Koebe transform $F$ of $f$) to obtain
\be\label{CP-extraeq5}
F(z)=\frac{f\left(\frac{z+\zeta}{1+\overline{\zeta} z}\right)-f(\zeta)}{h'(\zeta)(1-|\zeta |^2)}=:H(z)+\overline{G(z)}.
\ee
A calculation gives,
$$\frac{H''(0)}{2}=\frac{1}{2} \left \{(1-|\zeta|^{2})\frac{h''(\zeta)}{h'(\zeta)}-2\overline{\zeta}\right \}.
$$
Now, we consider the class ${\mathcal S}_{H_{2}}$ of all harmonic mappings
$f=h+\overline{g}\in{\mathcal S}_{H}$ satisfying
\be\label{CP-extraeq6}
\sup_{z\in\mathbb{D}}\left|(1-|z|^{2})\frac{h''(z)}{h'(z)}-2\overline{z}\right|<4 .
\ee
This inequality obviously holds if $h\in {\mathcal S}$ and $h$ is not the Koebe function $z/(1-e^{i\theta}z)^2$, $\theta\in \IR$.
Note that for the Koebe function the supremum turns out to be $4$. Our next two results are extension of \cite[Theorem 3]{Po}.

\begin{thm}\label{thm-1}
%For $K\geq1$, let $f\in{\mathcal S}_{H_{2}}$ be a $K$-quasiconformal
%harmonic mapping from $\mathbb{D}$ onto a bounded Pommerenke interior
%domain $G$.
Let $f\in{\mathcal S}_{H_{2}}(K, G)$, where $G$ is a bounded Pommerenke interior domain.
If there are positive constants $\delta_{2}\in(0,1)$
and $C_{5}$ such that, for each $\zeta\in\partial\mathbb{D}$
and for $0\leq \rho_{1}\leq\rho_{2}<1,$
\be\label{eqy}
\|D_{f}(\rho_{2}\zeta)\| \leq
C_{5}\left(\frac{1-\rho_{2}}{1-\rho_{1}}\right)^{\delta_{2}-1}\|D_{f}(\rho_{1}\zeta)\|,
\ee
then
$$\sup_{\zeta\in\mathbb{D}}\frac{1}{2\pi}\int_{\partial\mathbb{D}}\frac{\|D_{f}(\xi)\|}{\|D_{f}(\zeta)\|}
\frac{1-|\zeta|^{2}}{|\xi-\zeta|^{2}}|\,d\xi|<\infty.
$$
\end{thm}

We remark that if $K=1$, then Theorem \ref{thm-1} coincides with
\cite[Theorem 3]{Po}.

By using  similar reasoning as in the proof of Theorem \ref{thm-1}, one can easily get the following result which
replaces the assumption $f\in{\mathcal S}_{H_{2}}$ by a more general condition  $f\in{\mathcal S}_{H}$ and thus, we
omit its proof.

\begin{thm}\label{thm-1.02}
%For $K\geq1$, let $f\in{\mathcal S}_{H}$ be a $K$-quasiconformal
%harmonic mapping from $\mathbb{D}$ onto a bounded Pommerenke
%interior domain $G$.
Let $f\in{\mathcal S}_{H}(K, G)$, where $G$ is a bounded Pommerenke interior domain.
If there are  constants $C_{6}>0$, $C_{7}>0$,
$\delta_{3}>0$ and $\delta_{4}\in(0,1)$ such that, for each $\zeta\in\partial\mathbb{D}$
and for $0\leq \rho_{1}\leq\rho_{2}<1,$
$$C_{6}\left(\frac{1-\rho_{1}}{1-\rho_{2}}\right)^{\delta_{3}-1}\|D_{f}(\rho_{1}\zeta)\|\leq \|D_{f}(\rho_{2}\zeta)\|
\leq C_{7}\left(\frac{1-\rho_{2}}{1-\rho_{1}}\right)^{\delta_{4}-1}\|D_{f}(\rho_{1}\zeta)\|,
$$
then
$$\sup_{\zeta\in\mathbb{D}}\frac{1}{2\pi}\int_{\partial\mathbb{D}}\frac{\|D_{f}(\xi)\|}{\|D_{f}(\zeta)\|}
\frac{1-|\zeta|^{2}}{|\xi-\zeta|^{2}}|\,d\xi|<\infty.
$$
\end{thm}

Also, the following result easily follows from Theorem \ref{thm-1} and \cite[Theorem 1]{CP}.

\begin{cor}\label{cor-1}
For $K\geq1$, let $f\in{\mathcal S}_{H_{2}}\cap {\mathcal S}_{H}^{0}$ be a $K$-quasiconformal harmonic
mapping from $\mathbb{D}$ onto a bounded Pommerenke interior domain $G$. If $G$ is a radial John disk, then
$$\sup_{\zeta\in\mathbb{D}}\frac{1}{2\pi}\int_{\partial\mathbb{D}}\frac{\|D_{f}(\xi)\|}{\|D_{f}(\zeta)\|}
\frac{1-|\zeta|^{2}}{|\xi-\zeta|^{2}}|\,d\xi|<\infty.
$$
\end{cor}

The proofs of Theorems \ref{thm-CPx}-\ref{thm-1} will be presented in Section \ref{csw-sec2}.

\section{The proofs of the main results }\label{csw-sec2}
Let $\lambda_{\mathbb{D}}$ stand for the  {\it hyperbolic distance} (or {\it Poincar\'e distance}) on the unit disk $\ID$.
We have
$$\lambda_{\mathbb{D}}(z_{1}, z_{2})=\inf_{\gamma}\int_{\gamma}\frac{|dz|}{1-|z|^{2}}=
\tanh^{-1}\left|\frac{z_{1}-z_{2}}{1-\overline{z}_{1}z_{2}}\right|,
$$
where the infimum is taken over all smooth curves $\gamma$ in $\ID$ connecting $z_{1}\in \ID$ and $z_{2}\in\mathbb{D}$ (cf. \cite{Po1}).
% It is well-known that,
%for $z_{1}, z_{2}\in\mathbb{D}$,
%$$\lambda_{\mathbb{D}}(z_{1}, z_{2})=\frac{1}{2}\log\frac{1+|z_{1}-z_{2}|/|1-\overline{z}_{1}z_{2}|}{1-|z_{1}-z_{2}|/|1-\overline{z}_{1}z_{2}|},
%$$
%which is equivalent to
%$$\left|\frac{z_{1}-z_{2}}{1-\overline{z}_{1}z_{2}}\right|=\frac{e^{2\lambda_{\mathbb{D}}(z_{1},
%z_{2})}-1}{e^{2\lambda_{\mathbb{D}}(z_{1}, z_{2})}+1}=\tanh\lambda_{\mathbb{D}}(z_{1}, z_{2}).
%$$
In \cite{Small}, Sheil-Small proved that if
$f=h+\overline{g}\in{\mathcal S}_{H}$, then
\be\label{eqLem-B}
\frac{(1-|z|)^{\alpha-1}}{(1+|z|)^{\alpha+1}}\leq|h'(z)|\leq\frac{(1+|z|)^{\alpha-1}}{(1-|z|)^{\alpha+1}}
\ee
and
$$ \alpha :=\sup_{f\in{\mathcal S}_{H}}\frac{|h''(0)|}{2}<\infty.
$$
Unless otherwise stated, the number $\alpha$ will be used throughout the
discussion and is indeed called the order of the linear invariant
family ${\mathcal S}_{H}$ (see \cite{Small}).

%\begin{Lem}\label{Lem-B}
%Let $f=h+\overline{g}\in{\mathcal S}_{H}$ and $\alpha :=\sup_{f\in{\mathcal S}_{H}}\frac{|h''(0)|}{2}<\infty$. Then
%$$\frac{(1-|z|)^{\alpha-1}}{(1+|z|)^{\alpha+1}}\leq|h'(z)|\leq\frac{(1+|z|)^{\alpha-1}}{(1-|z|)^{\alpha+1}}.
%$$
%%where $\alpha$ is the same as in Theorem {\rm \ref{thm-c1}.}
%\end{Lem}

\begin{lem}\label{lem-1.1}
%For $K\geq1$, suppose that $f=h+\overline{g}\in{\mathcal S}_{H}$ is a $K$-quasiconformal harmonic mapping.
Suppose that $f\in{\mathcal S}_{H}(K)$. Then, for $z_{0},
z_{1}\in\mathbb{D}$,
$$\frac{1}{\alpha(1+K)}\left[1-e^{-2\alpha\lambda_{\mathbb{D}}(z_{1},z_{0})}\right]\leq
\frac{|f(z_{1})-f(z_{0})|}{(1-|z_{0}|^{2})|f_{z}(z_{0})|}
\leq\frac{K}{\alpha(1+K)}\left[e^{2\alpha\lambda_{\mathbb{D}}(z_{1},z_{0})}-1\right].
$$
In particular,
\be\label{cor-CP}
\frac{1}{\alpha(1+K)}\left[1-e^{-2\alpha\lambda_{\mathbb{D}}(z,0)}\right]\leq|f(z)|\leq
\frac{K}{\alpha(1+K)}\left[e^{2\alpha\lambda_{\mathbb{D}}(z,0)}-1\right], \quad z\in\ID.
\ee
%where $\alpha$ is the same as in Theorem {\rm \ref{thm-c1}.}
\end{lem}
\bpf
%Let $f=h+\overline{g}\in{\mathcal S}_{H}$ be a $K$-quasiconformal harmonic
%mapping, where $h$ and $g$ are analytic in $\mathbb{D}$.
By assumption $f=h+\overline{g}\in{\mathcal S}_{H}$ is a $K$-quasiconformal harmonic mapping, where $h$ and
$g$ are analytic in $\mathbb{D}$. Thus, by \eqref{eqLem-B}, we have
\be\label{CP-extraeq7}
\|D_{f}(z)\| \leq \frac{2K}{K+1}|h'(z)| \leq  \frac{2K}{K+1}\frac{(1+|z|)^{\alpha-1}}{(1-|z|)^{\alpha+1}}
\ee
and thus,  for $z\in\mathbb{D},$ we obtain
\beq\label{eq-1p}
\nonumber |f(z)|&\leq&\int_{[0,z]}\|D_{f}(\zeta)\|\, |d\zeta|\\
%\nonumber &\leq&\frac{2K}{K+1}\int_{[0,z]}|h'(\zeta)|\,|d\zeta|\\
&\leq&\frac{2K}{K+1}\int_{0}^{|z|}\frac{(1+\rho)^{\alpha-1}}{(1-\rho)^{\alpha+1}}\,d\rho
 =\frac{K}{\alpha(K+1)}\left[\Bigg(\frac{1+|z|}{1-|z|}\Bigg)^{\alpha}-1\right].
\eeq
%\beq\label{eq-1p}
%\nonumber |f(z)|&\leq&\int_{[0,z]}\|D_{f}(\zeta)\| |d\zeta|\\
%\nonumber &\leq&\frac{2K}{K+1}\int_{[0,z]}|h'(\zeta)|\,|d\zeta|\\
%\nonumber &\leq&\frac{2K}{K+1}\int_{0}^{|z|}\frac{(1+\rho)^{\alpha-1}}{(1-\rho)^{\alpha+1}}d\rho\\
%&=&\frac{K}{\alpha(K+1)}\left[\Bigg(\frac{1+|z|}{1-|z|}\Bigg)^{\alpha}-1\right].
%\eeq
%where $[0,z]$ is the radial line segment from $0$ to $z$.
On the other hand, let $\Gamma$ be the preimage under $f$ of the radial segment from $0$ to $f(z)$. Again, because
$$l(D_{f}(z))\geq\frac{2}{K+1}|h'(z)| \geq \frac{2}{K+1} \frac{(1-|z|)^{\alpha-1}}{(1+|z|)^{\alpha+1}},
$$
it follows  that
\be\label{eq-2p} |f(z)| \geq
\int_{\Gamma}l(D_{f}(\zeta))\,|d\zeta| \geq \frac{1}{\alpha
(K+1)}\left[1-\Bigg(\frac{1-|z|}{1+|z|}\Bigg)^{\alpha}\right].
\ee
%\beq\label{eq-2p}
%\nonumber |f(z)|&=&\int_{\Gamma}\left|f_{\zeta}(\zeta)d\zeta+f_{\overline{\zeta}}(\zeta)d\overline{\zeta}\right|\\
%\nonumber &\geq&\int_{\Gamma}l(D_{f}(\zeta))|d\zeta|\\
%%\nonumber &\geq&\frac{2}{K+1}\int_{\Gamma}|h'(\zeta)|\,|d\zeta|\\
%\nonumber
%&\geq&\frac{2}{K+1}\int_{0}^{|z|}\frac{(1-\rho)^{\alpha-1}}{(1+\rho)^{\alpha+1}}d\rho\\
%\nonumber &=&\frac{1}{\alpha
%(K+1)}\left[1-\Bigg(\frac{1-|z|}{1+|z|}\Bigg)^{\alpha}\right].
%\eeq

Let $z=\frac{z_{1}-z_{0}}{1-\overline{z}_{0}z_{1}}$ so that
$z_{1}=\frac{z+z_{0}}{1+\overline{z}_{0}z}$, where
$z_{0},z_{1}\in\mathbb{D}.$ Then, by assumption,
$$F(z)= \frac{f(z)-f(z_{0})}{(1-|z_{0}|^{2})h'(z_{0})}\in{\mathcal S}_{H}
$$
and is a $K$-quasiconformal harmonic mapping, i.e. $F\in{\mathcal S}_{H}(K)$. Applying (\ref{eq-1p}) and (\ref{eq-2p}) to $F$ gives us the desired
result if we take into account of the fact that
$$\left|\frac{z_{1}-z_{2}}{1-\overline{z}_{1}z_{2}}\right|=\frac{e^{2\lambda_{\mathbb{D}}(z_{1},
z_{2})}-1}{e^{2\lambda_{\mathbb{D}}(z_{1}, z_{2})}+1}=\tanh\lambda_{\mathbb{D}}(z_{1}, z_{2}).
$$
% we get
%$$\frac{1}{\alpha(1+K)}\left[1-e^{-2\alpha\lambda_{\mathbb{D}}(z_{1},z_{0})}\right]
%\leq|F(z)| \leq \frac{K}{\alpha(1+K)}\left[e^{2\alpha\lambda_{\mathbb{D}}(z_{1},z_{0})}-1\right],
%$$
%which implies that
%$$\frac{1}{\alpha(1+K)}\left[1-e^{-2\alpha\lambda_{\mathbb{D}}(z_{1},z_{0})}\right]\leq\frac{|f(z_{1})-f(z_{0})|}{(1-|z_{0}|^{2})|h'(z_{0})|}\leq
% \frac{K}{\alpha(1+K)}\left[e^{2\alpha\lambda_{\mathbb{D}}(z_{1},z_{0})}-1\right].
%$$
The proof of the lemma is complete.
\epf

%The following result easily follows from Lemma \ref{lem-1.1}.
%
%\begin{cor}\label{cor-CP}
%For $K\geq1$, suppose that $f\in{\mathcal S}_{H}$ is a $K$-quasiconformal harmonic  mapping.
%Then, for $z\in\mathbb{D}$,
%$$\frac{1}{\alpha(1+K)}\left[1-e^{-2\alpha\lambda_{\mathbb{D}}(z,0)}\right]\leq|f(z)|\leq
%\frac{K}{\alpha(1+K)}\left[e^{2\alpha\lambda_{\mathbb{D}}(z,0)}-1\right],
%$$
%where $\alpha$ is the same as in Theorem {\rm \ref{thm-c1}.}
%\end{cor}

\begin{lem}\label{lem-cp-u1}
Assume that $f\in{\mathcal S}_{H}(K)$. Then
$$\|D_{f}(z)\|\,|z|\leq\frac{C_{8}|f(z)|}{1-|z|} ~\mbox{ for $z\in\mathbb{D}$},
$$
where
\be\label{CP-extraeq4} C_{8}=2\alpha
K\sup_{z\in\mathbb{D}}\left\{\frac{|z|(1+|z|)^{\alpha-1}}{\left[(1+|z|)^{\alpha}-(1-|z|)^{\alpha}\right]}\right\} \geq K.
\ee
%and  $\alpha$ is the same as in Theorem {\rm \ref{thm-c1}.}
\end{lem}
\bpf  Suppose that $f=h+\overline{g}\in{\mathcal S}_{H}(K)$, where $h$ and $g$ are
analytic in $\mathbb{D}$. Next, for fixed $\zeta\in\mathbb{D}$, consider the Koebe transform $F$ of $f$ given by
\eqref{CP-extraeq5}.
%\be\label{eq-cp-u2}
%F(z)=\frac{f\left(\frac{z+\zeta}{1+\overline{\zeta}z}\right)-f(\zeta)}{(1-|\zeta|^{2})h'(\zeta)},~z\in\mathbb{D}.
%\ee
By assumption,  $F\in{\mathcal S}_{H}$ and is also a $K$-quasiconformal harmonic mapping.
By letting $z=-\zeta$ in \eqref{CP-extraeq5} and applying \eqref{cor-CP} to $F$, we obtain (since $f(0)=0$)
$$|F(-\zeta)|=\frac{|f(\zeta)|}{(1-|\zeta|^{2})|h'(\zeta)|}\geq\frac{1}{\alpha(1+K)}
\frac{\left[(1+|\zeta |)^{\alpha}-(1-|\zeta|)^{\alpha}\right]}{(1+|\zeta |)^{\alpha}}
$$
which gives
$$\frac {|h'(\zeta)|}{(1+K)|f(\zeta)|}\leq \frac{(1+|\zeta |)^{\alpha-1}}{\left[(1+|\zeta |)^{\alpha}-(1-|\zeta|)^{\alpha}\right]}\cdot \frac{\alpha }{1-|\zeta|}.
$$
Since this follows for each $\zeta\in\ID$, by the first inequality in \eqref{CP-extraeq7},  we easily have
$$\frac{\|D_{f}(z)\|\,|z|}{|f(z)|}\leq\frac{2K}{K+1}\frac{|h'(z)|\,|z|}{|f(z)|}\leq \frac{C_{8}}{1-|z|},
$$
where $C_8$ is given by \eqref{CP-extraeq4}.
%$$C_{8}=2\alpha K\sup_{z\in\mathbb{D}}\left\{\frac{|z|(1+|z|)^{\alpha-1}}{(1+|z|)^{\alpha}-(1-|z|)^{\alpha}}\right\}\geq K.
%$$
\epf

\begin{lem}\label{lem-cp-u3}
%For $K\geq1$, suppose that $f\in{\mathcal S}_{H}$ is a $K$-quasiconformal harmonic mapping.
Let $f\in{\mathcal S}_{H}(K)$ and, for any fixed
$\theta\in[0,2\pi]$, set
$$m_{f}(r,\theta)=\max_{\rho\in[0,r]}|f(\rho e^{i\theta})|,
$$
where $r\in[0,1)$.  Then, for $0<\rho_{0}\leq r<1$
and $0\leq\rho\leq r$, there is a constant $C_{9}>0$ which depends only on $\rho_{0}$ such that
\be\label{eq-cp-u}
\frac{|f(\rho e^{i\theta})|}{\rho}\leq C_{9}m_{f}(r,\theta),
\ee
where $\rho_{0}$ is a constant.
\end{lem}
\bpf Without loss of generality, we assume that $\theta=0.$ Clearly, \eqref{cor-CP} yields that
%\be\label{eq-cp-u4}
$$\lim_{\rho\rightarrow0^{+}}\frac{|f(\rho)|}{\rho}\leq\frac{K}{\alpha(1+K)}
\lim_{\rho\rightarrow0^{+}}\frac{e^{2\alpha\lambda_{\mathbb{D}}(\rho,0)}-1}{\rho}=\frac{2K}{1+K},
$$
%\ee
which implies that $f(\rho)/\rho$ is bounded in $[0,\rho_{0}]$,
where $\rho_{0}$ is a constant such that $0<\rho_{0}\leq r<1$.
%and $\alpha$ is the same as in Theorem {\rm \ref{thm-c1}.}
Hence there is a constant $C_{10}>0$ such that
\be\label{eq-cp-u6}
\frac{f(\rho)}{\rho}\leq C_{10}m_{f}(r,0)~\mbox{for}~\rho\in[0,\rho_{0}],
\ee
where $r\in[\rho_{0},1)$.  For $r\in[0,1)$, let
$$T(r)=\frac{K}{\alpha(1+K)}\left[1-e^{-2\alpha\lambda_{\mathbb{D}}(r,0)}\right].
$$
Then $T$ is increasing in $[0,1)$, which, together with \eqref{cor-CP}, yields that
\be\label{eq-cp-u5}
0<T(\rho_{0})\leq T(\rho)\leq f(\rho)\leq\frac{f(\rho)}{\rho} \leq
\frac{f(\rho)}{\rho_{0}}\leq\frac{1}{\rho_{0}}m_{f}(r,0)~\mbox{for}~
\rho\in[\rho_{0}, r],
\ee
where $r\in[\rho_{0},1)$. Therefore, (\ref{eq-cp-u}) follows from (\ref{eq-cp-u6}) and (\ref{eq-cp-u5}).
\epf

\begin{figure}[!ht]
\centering
\includegraphics{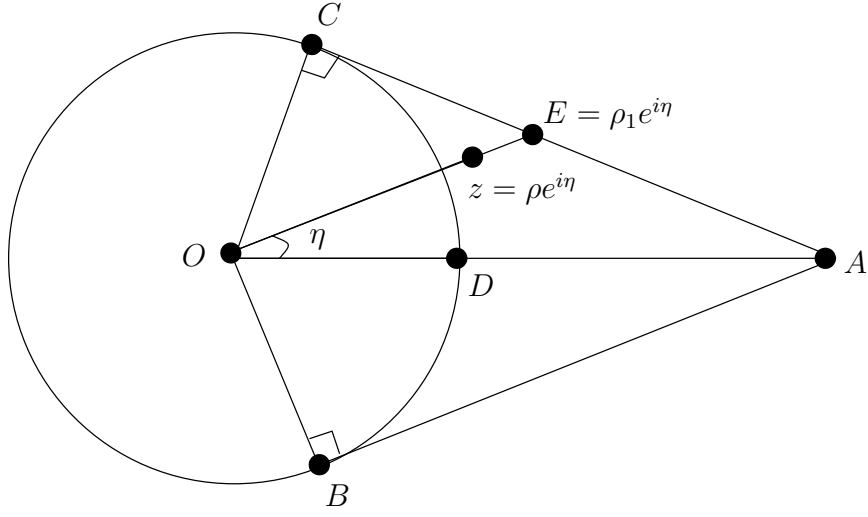}
\caption{Stolz-type domain} \label{fig1}
\end{figure}

\begin{lem}\label{lem-yp}
For $r\in(0,1)$, let $\Omega_{r}$ be the Stolz-type domain
consisting of the interior of the convex hull of the point $r$ and
the disk $\mathbb{D}_{r/4}$. Then,  for   $z=\rho
e^{i\eta}\in\Omega_{r}\backslash\mathbb{D}_{r/4},$
$$|\eta|\leq\frac{4\pi}{r\sqrt{15}}(r-\rho)<\frac{4\pi}{r\sqrt{15}}(1-\rho).
$$
\end{lem}
\bpf Assume without loss of generality that $\eta\geq0.$ Let $A$,
$D$, and $E$ represent the points $r$, $r/4$, and
$\rho_{1}e^{i\eta}$ (see Figure \ref{fig1}), respectively.
As $\angle OCA =\pi/2$, it is clear that
$$\sin\angle COA= \frac{\sqrt{15}}{4}, ~\cos\angle COA=\frac{1}{4},~
\sin\angle COE =\frac{\sqrt{\rho_{1}^{2}-\frac{r^{2}}{16}}}{\rho_{1}} ~\mbox{ and } \cos\angle COE=\frac{r}{4\rho_{1}}.
$$
Then, because $|(\sin \eta)/\eta| \geq 2/\pi$ for $|\eta|<\pi/2$, it follows that for $\eta \ge 0$
\begin{eqnarray*}
\frac{2\eta}{\pi}\leq\sin \eta&=&\sin\left(\angle COA-\angle
COE\right)\\&=&\sin\angle COA\cos\angle COE-\cos\angle COA\sin\angle
COE\\&=&\frac{\sqrt{15r^{2}}-\sqrt{16\rho_{1}^{2}-r^{2}}}{16\rho_{1}}\\
&=&\frac{r^{2}-\rho_{1}^{2}}{\rho_{1}\left(\sqrt{15}r+\sqrt{16\rho_{1}^{2}-r^{2}}\right)}
\end{eqnarray*}
Note that $\frac{r}{4}<\rho <\rho_1<r$ and, because
$$\frac{r^{2}-\rho_{1}^{2}}{\rho_{1}}<\frac{4}{r}\left (r^{2}-\rho_{1}^{2} \right ) <8(r-\rho_{1})  <8(r-\rho),
$$
the last above relation clearly implies that
$$\frac{2\eta}{\pi} <\frac{8(r-\rho)}{r\sqrt{15}}
$$
which gives the desired conclusion. Observe that $\frac{8}{r\sqrt{15}}(r-\rho)$ is less than $6/\sqrt{15}$ from which we also deduce that
$|\eta|< 3\pi/\sqrt{15}$.
%$$\eta\leq\frac{4\pi}{r\sqrt{15}}(r-\rho)\leq\frac{4\pi}{r\sqrt{15}}(1-\rho).
%$$
\epf

\subsection*{Proof of Theorem \ref{thm-CPx}}
Assume without loss of generality that $\theta=0$. For $r\in(0,1)$, we use
$\Omega_{r}$ to denote the Stolz-type domain, where $\Omega_{r}$ is
same as in Lemma \ref{lem-yp}. Let $z=\rho
e^{i\eta}\in\Omega_{r}\backslash\mathbb{D}_{r/4}$. Then, by Lemma
\ref{lem-yp}, there is a constant $C_{11}>0$ which depends only on $r$
such that
\be\label{eq-cp-u7}
|\eta|<C_{11}(1-\rho).
\ee
Suppose that $f=h+\overline{g}\in {\mathcal S}_{H}(K)$. By calculations, we get
\be\label{eq-cp-u8}
\log\frac{f(\rho e^{i\eta})}{\rho e^{i\eta}}-\log\frac{f(\rho)}{\rho}=i\int_{0}^{\eta}\left(\frac{\rho e^{it}h'(\rho e^{it})-
\rho e^{-it}\overline{g'(\rho e^{it})}}{f(\rho e^{it})}-1\right)dt.
\ee
Taking real part of \eqref{eq-cp-u8} on both sides, and then using (\ref{eq-cp-u7}), (\ref{eq-cp-u8}) and Lemma \ref{lem-cp-u1},
we see that there is a constant $C_{12}$ such that
\beq
\nonumber\log\frac{|f(\rho
e^{i\eta})|}{\rho}-\log\frac{|f(\rho)|}{\rho}&\leq&\int_{0}^{\eta}\frac{\rho\|D_{f}(\rho
e^{it})\|}{|f(\rho e^{it})|}\,dt\\ \nonumber
&\leq&C_{12}\int_{0}^{\eta}\frac{dt}{1-\rho}\leq C_{11}C_{12},
\eeq
which gives that
\be\label{eq-cp-u9}
|f(z)|=|f(\rho e^{i\eta})|\leq e^{C_{11}C_{12}}|f(\rho)|.
\ee
By \eqref{cor-CP}, we see that (\ref{eq-cp-u9}) also holds for $z\in\mathbb{D}_{r/4}$.
Then, by (\ref{eq-cp-u9}), there is constant $C_{13}$ such that
$f(\Omega_{r})$ is contained in $\mathbb{D}_{C_{13}m_{f}(r,0)}$,
which yields that
\be\label{eq-cp-u10}
\int_{\Omega_{r}}J_{f}(\zeta)\,dA(\zeta)\leq C_{13}^{2}m_{f}^{2}(r,0),
\ee
where $\zeta=x+iy,$  $dA=dxdy/\pi$ and $m_{f}(r,\theta)$ is defined as in Lemma \ref{lem-cp-u3}.

By \cite[Theorem 2]{MZ}, there is a constant $C_{14}$ such that
\be\label{eq-cp-u12}
\int_{0}^{1}(1-\rho)|H'(\rho)|^{2}\,d\rho\leq C_{14}\int_{\Omega_{1}}|H'(z)|^{2}\,dA(z),
\ee
where $H(z)$ is analytic in $\mathbb{D}$. For $r\in(0,1)$, let $H(z)=h(rz)$ for $z\in\mathbb{D}$.
Then, by (\ref{eq-cp-u12}), we  obtain
$$ \int_{0}^{r}(r-\rho)|h'(\rho)|^{2}\,d\rho\leq C_{15}\int_{\Omega_{r}}|h'(z)|^{2}\,dA(z),
$$
which implies that
\beq%\label{eq-cp-u13}
\nonumber
\int_{0}^{r}(r-\rho)\|D_{f}(\rho)\|^{2}\,d\rho
&\leq & C_{15}K\int_{\Omega_{r}}J_{f}(z)\,dA(z),\\
\label{eq-cp-u14} &\leq& C_{13}^{2}C_{15}Km_{f}^{2}(r,0) ~\mbox{ (by (\ref{eq-cp-u10})),}
\eeq
where $C_{15}>0$ is a constant.
%By (\ref{eq-cp-u10}) and (\ref{eq-cp-u13}), we have
%\be\label{eq-cp-u14}
%\int_{0}^{r}(r-\rho)\|D_{f}(\rho)\|^{2}\,d\rho\leq C_{13}^{2}C_{15}Km_{f}^{2}(r,0).
%\ee

By Lemmas \ref{lem-cp-u1} and \ref{lem-cp-u3}, for $r\in(1/2,1)$ and $\rho\in[0,r]$,
there is a constant $C_{16}$ such that
\begin{eqnarray}\label{eq-cp-u11}
\int_{0}^{r}\|D_{f}(\rho)\|^{2}\,d\rho&\leq&C_{16}\int_{0}^{r}\frac{m_{f}^{2}(r,0)}{(1-\rho)^{2}}\,d\rho
%~\mbox{(by Lemmas \ref{lem-cp-u1} and \ref{lem-cp-u3})}\\ &=&
= C_{16}m_{f}^{2}(r,0)\frac{r}{1-r}.
\end{eqnarray}
%which yields that
%\be\label{eq-cp-u11}
%(1-r)\int_{0}^{r}\|D_{f}(\rho)\|^{2}d\rho\leq C_{16}m_{f}^{2}(r,0).
%\eq
Writing $1-\rho =(1-r)+(r-\rho)$ and then, applying
(\ref{eq-cp-u14}) and (\ref{eq-cp-u11}), it follows  that
\beq\label{eq-cp-u15} \int_{0}^{r}(1-\rho)\|D_{f}(\rho)\|^{2}\,d\rho
%&=&\int_{0}^{r}(1-r)\|D_{f}(\rho)\|^{2}d\rho+ \int_{0}^{r}(r-\rho)\|D_{f}(\rho)\|^{2}d\rho\\
&\leq& (C_{13}^{2}C_{15}K+C_{16})m_{f}^{2}(r,0).
\eeq
Therefore, by (\ref{eq-cp-u15}), we conclude that
\begin{eqnarray*}
\ell_{f}(0,r)&\leq&\int_{0}^{r}\|D_{f}(\rho)\|\,d\rho\\
&\leq&\left(\int_{0}^{r}(1-\rho)\|D_{f}(\rho)\|^{2}\,d\rho\right)^{\frac{1}{2}}\left(\int_{0}^{r}\frac{d\rho}{1-\rho}\right)^{\frac{1}{2}}
\\&\leq&(C_{13}^{2}C_{15}K+C_{16})^{1/2}m_{f}(r,0)\left(\log\frac{1}{1-r}\right)^{\frac{1}{2}}.
\end{eqnarray*}

Now we prove the sharpness part. For any  $\tau\in(0,1/2)$, by
\cite{Je, Ke}, there is a function $h_0\in\mathcal S$ such that,
\be\label{eqis}
\ell_{h_{0}}(0,r)>C_{17}\left(\log\frac{1}{1-r}\right)^{\tau}~\mbox{as}~r\rightarrow1^{-},
\ee
where $C_{17}$ is a positive constant. Finally, consider
$$f_0(z)=h_0(z)+\frac{K-1}{K+1}\overline{h_0(z)}, \quad z\in\mathbb{D},
$$
and observe that $f_0\in{\mathcal S}_{H}$ and is a $K$-quasiconformal harmonic mapping. Consequently,
%\begin{eqnarray*}
$$\ell_{f_{0}}(0,r)=\int_{0}^{r}\left|h_{0}'(\rho)+\frac{K-1}{K+1}\overline{h_{0}'(\rho)}\right|d\rho\\
\geq\frac{2}{K+1}\int_{0}^{r}|h_{0}'(\rho)|\, d\rho=\frac{2}{K+1}\ell_{h_{0}}(0,r),
$$
which, together with (\ref{eqis}), implies that
$$\ell_{f_{0}}(0,r)>\frac{2C_{17}}{K+1}\left(\log\frac{1}{1-r}\right)^{\tau}~\mbox{as}~r\rightarrow1^{-}.
$$
The proof of this theorem is complete. \hfill $\Box$

\begin{lem}\label{lem-JD}
Let $f\in{\mathcal S}_{H}^{0}$. Then,  for
$\xi\in\partial\mathbb{D}$ and $0\leq\rho\leq r<1$,
\be\label{eq-w2}
\frac{(1-\rho^{2})\|D_{f}(\rho\xi)\|}{(1-r^{2})\|D_{f}(r\xi)\|}\leq
e^{2\alpha\lambda_{\mathbb{D}}(\rho,r)}.
\ee
%where $\alpha$ is the same as in Theorem {\rm \ref{thm-c1}.}
\end{lem}

\bpf Let $f=h+\overline{g}\in{\mathcal S}_{H}^{0}$, where $h$ and
$g$ are analytic in $\mathbb{D}$. For every $\mu\in\mathbb{D}$,
consider the affine mapping
$$f_{\mu}=f+\mu\overline{f}=(h+\mu g)+\overline{(g+\mu h)}.
$$
Clearly, $f_{\mu}\in{\mathcal S}_{H}$.  For a fixed $\zeta\in\mathbb{D}$, we
consider the Koebe transform  $F_{\mu}$ of $f_{\mu}$ as given by \eqref{CP-extraeq5}. Then
we can write $F_{\mu}=H_\mu+\overline{G_\mu}$
%$$K(z)=\frac{f_{\mu}(\frac{z+\zeta}{1+\overline{\zeta}z})-f_{\mu}(\zeta)}{(1-|\zeta|^{2})(h'(\zeta)+\mu g'(\zeta))}
%=H(z)+\overline{G(z)},
%$$
which again belongs to ${\mathcal S}_{H}$ and obviously,
%By elementary calculations, we get
%$$H(z)=z+A_{2}(\zeta)z^{2}+A_{3}(\zeta)z^{3}+\cdots,
%$$
%where
$$\frac{H_\mu ''(0)}{2} =A_{2}(\zeta)=\frac{1}{2}(1-|\zeta|^{2})\frac{h''(\zeta)+\mu g''(\zeta)} {h'(\zeta)+\mu g'(\zeta)}-\overline{\zeta}.
$$
Since $|A_{2}|\leq \alpha,$ we see that
\beq\label{eqt-7} \nonumber
\left|\frac{\partial}{\partial\rho}\log\big[(1-\rho^{2})(h'(\rho\xi)+\mu
g'(\rho\xi))\big]\right| &=&\left|\frac{h''(\rho\xi)+\mu
g''(\rho\xi)}{h'(\rho\xi)+\mu
g'(\rho\xi)}-\frac{2\rho\overline{\xi}} {1-\rho^{2}}\right|\\
\nonumber & \leq&\frac{2\alpha}{1-\rho^{2}},
\eeq
where $\xi\in\partial\mathbb{D}$. Integration leads to
$$\frac{(1-r^{2})|h'(r\xi)+\mu g'(r\xi)|}{(1-\rho^{2})|h'(\rho\xi)+\mu g'(\rho\xi)|}
\geq \left(\frac{1-r}{1+r}\cdot\frac{1+\rho}{1-\rho}\right)^{\alpha},
$$
which gives
%\be\label{eq-w1}
$$(1-\rho^{2})|h'(\rho\xi)+\mu g'(\rho\xi)|\leq e^{2\alpha\lambda_{\mathbb{D}}(\rho,r)}(1-r^{2})|h'(r\xi)+\mu g'(r\xi)|
$$
and the desired inequality (\ref{eq-w2}) follows from this and the arbitrariness of $\mu$.
\epf

%In \cite{M-1},  Mateljevi\'c  proved that if $f\in{\mathcal S}_{H}^{0}(K)$, then
%%\be\label{eq-M}
%$$d_{\Omega}(f(z))\geq\frac{\|D_{f}(z)\|(1-|z|^{2})}{16K} ~\mbox{ for $z\in\mathbb{D},$}
%$$
%where $\Omega=f(\mathbb{D})$ (see also \cite{Mi,M-2}). We wish to show this indeed holds for
%$f\in{\mathcal S}_{H}(K)$.

We remark that Mateljevi\'c \cite{M-1} (see also \cite{Mi,M-2}) proved the following lemma for
$f\in{\mathcal S}_{H}^{0}(K)$ instead of $f\in{\mathcal S}_{H}(K)$. That is, the normalization condition on $f$, namely,
$f_{\overline{z}}(0)=0$, is not necessary.

\begin{lem}\label{lem-ch-1}
If $f\in{\mathcal S}_{H}(K)$ and $\Omega=f(\mathbb{D})$,  then
\be\label{eq-ch}
d_{\Omega}(f(z))\geq\frac{\|D_{f}(z)\|(1-|z|^{2})}{16K}  ~\mbox{ for $z\in\mathbb{D}$}.
\ee
\end{lem}
\bpf Let $f=h+\overline{g}\in{\mathcal S}_{H}(K)$, where $h$ and
$g$ are analytic in $\mathbb{D}$. Then the affine mapping $f_0$ defined by
$$f_{0}(z)=\frac{f(z)-\overline{g'(0)}\overline{f(z)}}{1-|g'(0)|^{2}}
$$
belongs to ${\mathcal S}_{H}^{0}$. By \cite[Theorem 4.4]{Clunie-Small-84}, we have
\be\label{eq-ch-1}
\frac{|f(z)|}{1-|g'(0)|}\geq|f_{0}(z)|=\frac{|f(z)-\overline{g'(0)f(z)}|}
{1-|g'(0)|^{2}}\geq \frac{|z|}{4(1+|z|)^{2}}, \quad z\in\ID.
\ee
Recall again, for any fixed $\zeta\in\mathbb{D},$  the Koebe transform $F$ of $f$ given by
\eqref{CP-extraeq5}
%$$F(z)=\frac{f\left(\frac{z+\zeta}{1+\overline{\zeta}z}\right)-f(\zeta)}{(1-|\zeta|^{2})h'(\zeta)},~z\in\mathbb{D}.
%$$
belongs to ${\mathcal S}_{H}$ and $F$ is again a $K$-quasiconformal harmonic mapping. As a result,  (\ref{eq-ch-1}) applied to $F$ shows that
\begin{eqnarray*}
\bigg|f\left(\frac{z+\zeta}{1+\overline{\zeta}z}\right)-f(\zeta)\bigg|
&\geq& (1-|\zeta|^{2})|h'(\zeta)|(1-|F_{\overline{z}}(0)|) \frac{|z|}{4(1+|z|)^{2}}\\
&\geq& (1-|\zeta|^{2})|h'(\zeta) \left (\frac{2}{K+1}\right ) \frac{|z|}{4(1+|z|)^{2}}\\
&\geq&\frac{(1-|\zeta|^{2})\|D_{f}(\zeta)\|}{K}\frac{|z|}{4(1+|z|)^{2}},
\end{eqnarray*}
which implies that
$$d_{\Omega}(f(\zeta))=\liminf_{|z|\rightarrow1^{-}}
\frac{\bigg|f\left(\frac{z+\zeta}{1+\overline{\zeta}z}\right)-f(\zeta)\bigg|}
{|z|}\geq\frac{\|D_{f}(\zeta)\|(1-|\zeta|^{2})}{16K}.
$$
The proof of this Lemma is complete. \epf

\begin{Lem}{\rm (\cite[Lemma 2]{CP})} \label{Lem-CP}
Let $a_{1}, a_{2}$ and $a_{3}$ be positive constants and let
$0<|z_{0}|=1-\delta_{5}$, where $\delta_{5}\in(0,1)$. If
$f\in{\mathcal S}_{H}$,
$0\leq1-a_{2}\delta_{5}\leq|z|\leq1-a_{1}\delta_{5}$ and $|\arg z-
\arg z_{0}|\leq a_{3}\delta_{5}$, then
$$\frac{1}{M(a_{1},a_{2},a_{3})}\|D_{f}(z_{0})\|\leq\|D_{f}(z)\|\leq
M(a_{1},a_{2},a_{3})\|D_{f}(z_{0})\|,
$$
where
$M(a_{1},a_{2},a_{3})=2e^{(1+\alpha)\left(a_{3}+\frac{1}{2}\log\frac{2a_{2}-a_{1}}{a_{1}}\right)}.$
%and  $\alpha$ is the same as in Theorem {\rm \ref{thm-c1}.}
\end{Lem}

\subsection*{Proof of Theorem \ref{thm-J}} Let $f\in{\mathcal S}_{H}^{0}(K)$. First we show that ${\rm (ii)}\Rightarrow {\rm (i)}$.   We assume that
\be\label{eq-w3}
\frac{(1-\rho^{2})\|D_{f}(\rho\zeta)\|}{(1-r^{2})\|D_{f}(r\zeta)\|}\leq\beta<1~\mbox{for}~
\rho=\frac{x+r}{1+xr},~|\zeta|=1,
\ee
uniformly on $r$ and $\zeta$.  Define $x_1=x$ and $x_k$ for $k=2,3,\ldots$, by
$$\frac{1+x_{k}}{1-x_{k}}=\left(\frac{1+x}{1-x}\right)^{k}, \mbox{ i.e., }~ x_{k+1}=\frac{x+x_{k}}{1+xx_{k}}.
$$
Note that $\rho >r$ and thus, $x_{k+1}>x_{k}$. Consequently, by (\ref{eq-w3}), we have
\be\label{eq-w4}
\frac{(1-x_{k+1}^{2})\|D_{f}(x_{k+1})\|} {(1-x_{k}^{2})\|D_{f}(x_{k})\|}\leq\beta<1.
\ee
Let $\delta_{6}\in(0,1)$ such that
$$\beta<\left(\frac{1-x}{1+x}\right)^{\delta_{6}}.
$$
Then, for  $j<k$, by (\ref{eq-w4}),
\beq\label{eq-w5}
\nonumber \frac{(1-x_{k}^{2})\|D_{f}(x_{k})\|}
{(1-x_{j}^{2})\|D_{f}(x_{j})\|}&=&\frac{(1-x_{k}^{2})\|D_{f}(x_{k})\|}
{(1-x_{k-1}^{2})\|D_{f}(x_{k-1})\|}\\
\nonumber&\times&\frac{(1-x_{k-1}^{2})\|D_{f}(x_{k-1})\|}
{(1-x_{k-2}^{2})\|D_{f}(x_{k-2})\|}
\times\cdots\times\frac{(1-x_{j+1}^{2})\|D_{f}(x_{j+1})\|}
{(1-x_{j}^{2})\|D_{f}(x_{j})\|}\\ \nonumber &\leq& \beta^{k-j}\\
\nonumber &<&\left(\frac{1-x_{k}}{1+x_{k}}\right)^{\delta_{6}}
\left(\frac{1-x_{j}}{1+x_{j}}\right)^{-\delta_{6}}\\
&\leq&\left(\frac{1-x_{k}}{1-x_{j}}\right)^{\delta_{6}}.
\eeq
By calculations, for $k=\{1,2,\ldots\}$,
$$\lambda_{\mathbb{D}}(x_{k},x_{k+1})=\lambda_{\mathbb{D}}(0,x),
$$
which, together with (\ref{eq-w5}) and Lemma \ref{lem-JD}, yields
that there is a constant $C_{18}>0$ such that
\be\label{eq-w6}
\frac{\|D_{f}(\rho\zeta)\|}{\|D_{f}(r\zeta)\|}\leq
C_{18}\left(\frac{1-\rho}{1-r}\right)^{\delta_{6}-1}.
\ee
Hence, by (\ref{eq-w6}) and \cite[Theorem 1]{CP}, we conclude that $\Omega$ is
a radial John disk.

${\rm (i)}\Rightarrow {\rm (ii)}$.   Suppose that $\Omega=f(\mathbb{D})$ is a
radial John disk. Then, by \cite[Theorem 1]{CP},  there are constants
$C_{19}>0$ and $\delta_{7}\in(0,1)$ such that, for each
$\zeta\in\partial\mathbb{D}$ and for $0\leq\rho\leq r<1,$
$$\frac{(1-\rho^{2})\|D_{f}(\rho\zeta)\|}{
(1-r^{2})\|D_{f}(r\zeta)\|}\leq C_{19}\left(\frac{1-\rho}{1-r}\right)^{\delta_{7}}
 =C_{19}\left(\frac{1-x}{1+rx}\right)^{\delta_{7}}\leq C_{19}(1-x)^{\delta_{7}}.
$$
It is not difficult to see that
$C_{19}(1-x)^{\delta_{7}}<1$ by taking $x$ sufficiently close to 1.

Next we show that ${\rm (i)}\Rightarrow {\rm (iii)}$. For $z=re^{i\theta}\in\mathbb{D}$ and $w=r_{1}e^{i\theta_{1}}\in B(z)$,
by \cite[Theorem 1]{CP} and Lemma \Ref{Lem-CP}, we see that there are positive constants
$C_{20}$,  $C_{21}$ and $\delta_{8}\in(0,1)$ such that
\begin{eqnarray*}
|f(z)-f(w)|&\leq&|f(re^{i\theta})-f(re^{i\theta_{1}})|+|f(re^{i\theta_{1}})-f(r_{1}e^{i\theta_{1}})|\\
&\leq&r\int_{\gamma'}\|D_{f}(re^{it})\|\,dt+\int_{r}^{r_{1}}\|D_{f}(\rho
e^{i\theta_{1}})\|\,d\rho\\
&\leq&C_{20}r\int_{\gamma'}\|D_{f}(re^{i\theta})\|\,dt
+C_{21}\int_{r}^{r_{1}}\|D_{f}(re^{i\theta})\| \left(\frac{1-\rho}{1-r}\right)^{\delta_{8}-1}\,d\rho\\
&\leq&C_{20}
r\|D_{f}(re^{i\theta})\|\ell(\gamma')+\frac{C_{21}}{\delta_{8}}\|D_{f}(re^{i\theta})\|(1-r)\\
&\leq&\left(2\pi
C_{20}+\frac{C_{21}}{\delta_{8}}\right)\|D_{f}(re^{i\theta})\|(1-r),
\end{eqnarray*}
which gives that
$$\sup_{z\in\mathbb{D},~w\in B(z)}\frac{|f(z)-f(w)|}{(1-|z|^{2})\|D_{f}(z)\|}<\infty,
$$
where $\gamma'$ is the smaller subarc of $\partial\mathbb{D}_{r}$ between
$re^{i\theta}$ and $re^{i\theta_{1}}$.

 Finally, we prove ${\rm (iii)}\Rightarrow {\rm (i)}$. For $z\in\mathbb{D}$ and $w_{1},w_{2}\in B(z)$,  there is a positive
constant $C_{22}$ such that
\begin{eqnarray*}
|f(w_{1})-f(w_{2})|&\leq&|f(w_{1})-f(z)|+|f(w_{2})-f(z)|\\&\leq&C_{22}(1-|z|^{2})\|D_{f}(z)\|\\
&\leq&16KC_{22}d_{\Omega}(f(z)) ~\mbox{ (by Lemma \ref{lem-ch-1}),}
\end{eqnarray*}
which implies that
\be\label{eqwr1}
\diam f(B(z))\leq 16KC_{22}d_{\Omega}(f(z)).
\ee
By (\ref{eqwr1}) and \cite[Theorem 2]{CP}, we conclude that $\Omega$
is a radial John disk. The proof of this theorem is complete.
\hfill $\Box$

%\begin{Thm} $($\cite[Theorem 3]{CPRW}$)$\label{ThmC}
%Let $f\in{\mathcal S}_H^{0}$. Then there is a positive constant
%$C_{18}<\infty$ such that for $\xi\in\partial\mathbb{D}$ and $0\leq
%r_{3}\leq r_{4}<1$,
%$$\|D_{f}(r_{4}\xi)\|\geq\frac{1}{2^{1+C_{18}}}\|D_{f}(r_{3}\xi)\|\left(\frac{1-r_{4}}{1-r_{3}}\right)^{C_{18}-1}.$$
%\end{Thm}

\subsection*{Proof of Theorem \ref{thm-c1}} Let $f=h+\overline{g}\in{\mathcal S}_H^{0}(K,\Omega)$, where $\Omega$ is a radial John disk.
Assume that $a_{1}=re^{i\theta}$ and
$r_{1}e^{i\theta_{1}},r_{2}e^{i\theta_{2}}\in B(a_{1})$ with
$r_{1}\leq r_{2}$, where $r=|a_{1}|.$ Since $\Omega$ is a radial
John disk $\Omega$, by  \cite[Theorem 1]{CP}, we see that there are
constants $C_{23}>0$ and $\delta_{9}\in(0,1)$ such that for each
$\zeta\in\partial\mathbb{D}$ and for $0\leq r\leq\rho<1,$
\be\label{eq-p-11}
\|D_{f}(\rho\zeta)\|\leq C_{23} \|D_{f}(r\zeta)\|\left(\frac{1-\rho}{1-r}\right)^{\delta_{9}-1}.
\ee
Then, by (\ref{eq-p-11}) and Lemma  \Ref{Lem-CP}, there is a positive constant $C_{24}$ such that
\begin{eqnarray*}
|f(r_{2}e^{i\theta_{2}})-f(r_{1}e^{i\theta_{1}})|
&\leq&|f(r_{2}e^{i\theta_{2}})-f(re^{i\theta_{2}})|+|f(r_{1}e^{i\theta_{1}})-f(re^{i\theta_{1}})|\\
&&+|f(re^{i\theta_{2}})-f(re^{i\theta_{1}})|\\
&\leq&\int_{r}^{r_{2}}\|D_{f}(\rho e^{i\theta_{2}})\|d\rho
+\int_{r}^{r_{1}}\|D_{f}(\rho e^{i\theta_{1}})\|\,d\rho +J\\
%&& +r\int_{\gamma_{0}}\|D_{f}(r e^{it})\|\,dt\\
&\leq& C_{23}\left [\int_{r}^{r_{2}}\|D_{f}(re^{i\theta})\|\left(\frac{1-\rho}{1-r}\right)^{\delta_{9}-1}\,d\rho \right .
\\&&
\left . + \int_{r}^{r_{1}}\|D_{f}(re^{i\theta})\|\left(\frac{1-\rho}{1-r}\right)^{\delta_{9}-1}\,d\rho\right ]+J\\
&\leq&\frac{2C_{23}}{\delta_{9}}\|D_{f}(re^{i\theta})\|(1-r) +J,
\end{eqnarray*}
where
\begin{eqnarray*}
J& =&r\int_{\gamma_{0}}\|D_{f}(r e^{it})\|\,dt\leq C_{24}r\int_{\gamma_{0}}\|D_{f}(r e^{i\theta})\|\,dt\\
&\leq& C_{24}|\theta_{2}-\theta_{1}|\|D_{f}(re^{i\theta})\|\\
&\leq& 2\pi C_{24}\|D_{f}(re^{i\theta})\|(1-r),
\end{eqnarray*}
where $\gamma_{0}$ is the smaller subarc of $\partial\mathbb{D}_{r}$ between $re^{i\theta_{1}}$ and $re^{i\theta_{2}}$.
Combining the last two inequalities shows that
$$|f(r_{2}e^{i\theta_{2}})-f(r_{1}e^{i\theta_{1}})| \leq \left(\frac{2C_{23}}{\delta_{9}}+2\pi C_{24}\right)\|D_{f}(re^{i\theta})\|(1-r)
$$
Hence there is a constant $C_{25}>0$  such that
\be\label{eq-p-12}
\diam B(a_{1})\leq C_{25}(1-|a_{1}|)\|D_{f}(a_{1})\|.
\ee
Moreover, by Lemmas \ref{lem-ch-1} and \Ref{Lem-CP}, we see that there is a constant $C_{26}>0$  such that
\beq\label{eq-p-13}
\nonumber\diam f(B(a_{2}))&\geq&d_{\Omega}(f(a_{2}))\\
\nonumber &\geq&\frac{1}{16K}(1-|a_{2}|^{2})\|D_{f}(a_{2})\|\\
\nonumber &\geq&\frac{1}{16K}(1-|a_{2}|)\|D_{f}(a_{2})\|\\
&\geq&\frac{C_{26}}{16K}(1-|a_{2}|)\|D_{f}(|a_{2}|e^{i\theta})\|.
\eeq
By (\ref{eq-p-12}), (\ref{eq-p-13}) and Lemma \ref{lem-JD}, we conclude that
\begin{eqnarray*}
\frac{\diam f(B(a_{1}))}{\diam
f(B(a_{2}))}&\leq&\frac{16KC_{25}}{C_{26}}\frac{(1-|a_{1}|)\|D_{f}(a_{1})\|}{(1-|a_{2}|)\|D_{f}(|a_{2}|e^{i\theta})\|}\\
&\leq&\frac{2^{5+\alpha}KC_{25}}{C_{26}}\frac{(1-|a_{1}|)}{(1-|a_{2}|)}\left(\frac{1-|a_{1}|}{1-|a_{2}|}\right)^{\alpha-1}\\
&=&\frac{2^{5+\alpha}KC_{25}}{C_{26}}\left(\frac{1-|a_{1}|}{1-|a_{2}|}\right)^{\alpha}
\end{eqnarray*} %where $C_{18}$ is the same as in Theorem \Ref{ThmC}.
and the proof of the theorem is complete. \hfill $\Box$

\vspace{6pt}

\begin{lem}\label{lem-1.2}
For $K\geq1$, suppose that $f\in{\mathcal S}_{H}(K)$.
Let $a_{1}, a_{2}$ and $a_{3}$ be positive constants and let $0<|z_{0}|=1-\delta$, where
$\delta\in(0,1)$. Suppose further that $0\leq1-a_{2}\delta\leq|z|\leq1-a_{1}\delta$ and
$|\arg z- \arg z_{0}|\leq a_{3}\delta$. Then
$$|f(z)-f(z_{0})|\leq\frac{K}{\alpha(1+K)}
\left[\Bigg(\frac{M(a_{1},a_{2},a_{3})}{2}\Bigg)^{\frac{2\alpha}{1+\alpha}}-1\right](1-|z_{0}|^{2})|f_{z}(z_{0})|,
$$
where $M(a_{1},a_{2},a_{3})$ is defined in Lemma \Ref{Lem-CP}.
%$\alpha$ is the same as in Theorem {\rm \ref{thm-c1}.}
\end{lem}
\bpf Follows from \cite[Lemma 2]{CP}, but for the sake of completeness, we include certain details.

\begin{figure}[!ht]
\centering
\includegraphics{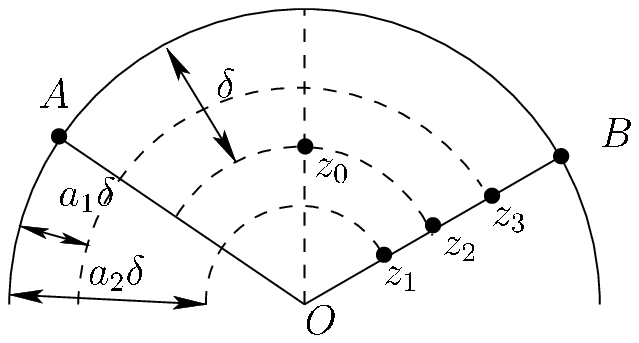}
\caption{ } \label{fig2}
\end{figure}

Let $\angle AOB=2a_{3}\delta$ and $z_{1}, z_{2}, z_{3}$ line in
the line $OB$ with $|z_{1}|\leq|z_{2}|=|z_{0}|\leq|z_{3}|$ (see
Figure \ref{fig2}). Clearly the distance from $O$ to $B$ is less than $1$. %(cf. \cite[Lemma 2]{CP}).
Then the length of the circular arc from $z_{0}$ to $z_{2}$ is less than $a_{3}\delta$.
As in \cite[Lemma 2]{CP}, it  is easy to see that
$$\lambda_{\mathbb{D}}(z_{0}, z_{2})<a_{3}, ~ \left|\frac{z_{3}-z_{1}}{1-\overline{z}_{1}z_{3}}\right| \leq\frac{a_{2}-a_{1}}{a_{2}}
~\mbox{ and }~
\lambda_{\mathbb{D}}(z_{0}, z) \leq a_{3}+\frac{1}{2}\log\frac{2a_{2}-a_{1}}{a_{1}}.
$$
%%<\frac{a_{3}\delta}{1-(1-\delta)^{2}}=\frac{a_{3}}{2-\delta}<a_{3}
%$$
%and
%$$\left|\frac{z_{3}-z_{1}}{1-\overline{z}_{1}z_{3}}\right|=\frac{1-a_{1}\delta-(1-a_{2}\delta)}
%{1-(1-a_{1}\delta)(1-a_{2}\delta)}=\frac{a_{2}-a_{1}}{a_{2}+a_{1}(1-a_{2}\delta)}\leq\frac{a_{2}-a_{1}}{a_{2}}.
%$$
%so that
%$$
%\begin{eqnarray*}
%\lambda_{\mathbb{D}}(z_{0}, z)&\leq&\lambda_{\mathbb{D}}(z_{0},
%z_{2})+\lambda_{\mathbb{D}}(z_{2}, z_{1})\\
%&\leq&\lambda_{\mathbb{D}}(z_{0}, z_{2})+\lambda_{\mathbb{D}}(z_{1},
%z_{3})\\
%&\leq&a_{3}+\frac{1}{2}\log\frac{2a_{2}-a_{1}}{a_{1}}.
%\end{eqnarray*}
The desired conclusion follows if we use Lemma \ref{lem-1.1}.
%, we see that
%$$|f(z)-f(z_{0})|\leq\frac{K}{\alpha(1+K)}
%\left[\Bigg(\frac{M(a_{1},a_{2},a_{3})}{2}\Bigg)^{\frac{2\alpha}{1+\alpha}}-1\right](1-|z_{0}|^{2})|f_{z}(z_{0})|,
%$$
%where $M(a_{1},a_{2},a_{3})$ is defined as in the statement.
%The proof of this lemma is complete.
\epf

The following result is an easy consequence of Lemmas \ref{lem-ch-1} and \ref{lem-1.2}.

\begin{cor}
%For $K\geq1$, suppose that $f\in{\mathcal S}_{H}$ is a
%$K$-quasiconformal harmonic mapping. Let $a_{1}, a_{2}$ and $a_{3}$
%be positive constants and let $0<|z_{0}|=1-\delta$, where
%$\delta\in(0,1)$. If $0\leq1-a_{2}\delta\leq|z|\leq1-a_{1}\delta$
%and $|\arg z- \arg z_{0}|\leq a_{3}\delta$, then
Under the hypotheses of Lemma \ref{lem-1.2}, we also have
$$|f(z)-f(z_{0})|\leq\frac{16K^{2}}{\alpha(1+K)}
\left[\Bigg(\frac{M(a_{1},a_{2},a_{3})}{2}\Bigg)^{\frac{2\alpha}{1+\alpha}}-1\right]d_{f(\mathbb{D})}(f(z_{0})),
$$
where $M(a_{1},a_{2},a_{3})$ is defined in Lemma \Ref{Lem-CP}.% and
%$\alpha$ is the same as in Theorem {\rm \ref{thm-c1}.}
\end{cor}

\subsection*{Proof of Theorem \ref{thm-c3}} Let $z=re^{i\theta}$,
$\sigma=|\zeta_{1}-\zeta_{2}|$ and  $\zeta_{j}=r_{j}e^{i\theta_{j}}$ ($j=1,2$) with $r_{1}\leq r_{2}$. \\
{\rm $\mathbf{Step~ 1.} $} If $r\leq\rho=1-2\sigma\leq r_{1}\leq r_{2}$, then
\begin{eqnarray*}
|\zeta_{1}-\zeta_{2}|&=&|r_{1}e^{i\theta_{1}}-r_{2}e^{i\theta_{2}}|\\
&=&\sqrt{r_{1}^{2}+r_{2}^{2}-2r_{1}r_{2}\cos(\theta_{1}-\theta_{2})}\\
&=&\sqrt{(r_{1}-r_{2})^{2}+4r_{1}r_{2}\sin^{2}\left(\frac{\theta_{1}-\theta_{2}}{2}\right)}\\
&\geq&2\sqrt{r_{1}r_{2}}\left|\sin\frac{\theta_{1}-\theta_{2}}{2}\right|\\
&\geq&\frac{2\rho|\theta_{1}-\theta_{2}|}{\pi},
\end{eqnarray*}
which, together  with \cite[Theorem 2]{CP}, Lemmas \ref{lem-ch-1}
and \ref{lem-1.2}, imply that there are positive constants
$C_{27}$, $C_{28}$, $C_{29}$, $C_{30}$, and $\delta_{10}\in(0,1)$ such that
\begin{eqnarray*}
|f(\zeta_{1})-f(\zeta_{2})|&\leq&|f(\zeta_{1})-f(\rho
e^{i\theta_{1}})|+|f(\zeta_{2})-f(\rho e^{i\theta_{2}})|+|f(\rho
e^{i\theta_{1}})-f(\rho e^{i\theta_{2}})|\\
&\leq&C_{27}\big[(1-\rho)\|D_{f}(\rho
e^{i\theta_{1}})\|+(1-\rho)\|D_{f}(\rho e^{i\theta_{2}})\|\big]\\
&&+\rho\int_{\gamma_{1}}\|D_{f}(\rho e^{it})\|\,dt~ \mbox{ (by Lemma \ref{lem-1.2}})\\
&\leq&C_{27}\big[(1-\rho)\|D_{f}(\rho
e^{i\theta_{1}})\|+(1-\rho)\|D_{f}(\rho e^{i\theta_{2}})\|\big]\\
&&+\rho\int_{\gamma_{1}}\|D_{f}(z)\|\left(\frac{1-\rho}{1-|z|}\right)^{\delta_{10}-1}\,dt~ \mbox{ (by \cite[Theorem 2]{CP}})\\
&\leq&\|D_{f}(z)\|\left(\frac{1-\rho}{1-|z|}\right)^{\delta_{10}-1}\big[C_{28}(1-\rho)+
C_{29}\ell(\gamma_{1})\big]\\
%&=&\|D_{f}(z)\|\left(\frac{1-\rho}{1-|z|}\right)^{\delta_{10}-1}\big[C_{28}(1-\rho)+
%C_{29}\rho|\theta_{1}-\theta_{2}|\big]\\
&\leq&\|D_{f}(z)\|\left(\frac{1-\rho}{1-|z|}\right)^{\delta_{10}-1}\big[C_{28}(1-\rho)+
\frac{C_{29}\pi\sigma}{2}\big]\\
&\leq&
C_{30}d_{f}(z)\left(\frac{|\zeta_{1}-\zeta_{2}|}{1-|z|}\right)^{\delta_{10}}~(\mbox{ by
Lemma \ref{lem-ch-1}})
\end{eqnarray*}
where $\gamma_{1}$ is the smaller subarc of $\partial\mathbb{D}_{\rho}$ between $\rho e^{i\theta_{1}}$ and $\rho
e^{i\theta_{2}}$, and
$$\ell(\gamma_{1}) =\rho|\theta_{1}-\theta_{2}|\leq \frac{\pi\sigma}{2} .
$$

\noindent
{\rm $\mathbf{Step~ 2.} $} If $r_{1}<\rho=1-2\sigma$, then, by Lemma
\ref{lem-1.2},  there are positive constants $C_{31}$ and $C_{32}$
such that
\be\label{eq-p14}
\|D_{f}(\zeta)\|\leq C_{31}\|D_{f}(\zeta_{1})\|\leq
C_{32}\|D_{f}(\rho e^{i\theta_{1}})\|,
\ee
where $|\zeta-\zeta_{1}|\leq\sigma$.  %By (\ref{eq-p14}), \cite[Theorem 2]{CP} and Lemma \ref{lem-ch-1},
We see that there are positive constants $C_{33}$ and  $\delta_{11}\in(0,1)$ such that
\begin{eqnarray*}
|f(\zeta_{1})-f(\zeta_{2})|&\leq&\int_{[\zeta_{1},\zeta_{2}]}\|D_{f}(\zeta)\|\,|d\zeta|\\
&\leq&C_{32}\|D_{f}(\rho e^{i\theta_{1}})\|\,|\zeta_{1}-\zeta_{2}|~\mbox{ (by (\ref{eq-p14}) })\\
&\leq&C_{33}\|D_{f}(z)\|\,|\zeta_{1}-\zeta_{2}|\left(\frac{1-\rho}{1-r}\right)^{\delta_{11}-1}~\mbox{ (\cite[Theorem 2]{CP})}\\
&\leq&2^{3+\delta_{11}}KC_{33}d_{\Omega}(f(z))\left(\frac{|\zeta_{1}-\zeta_{2}|}{1-|z|}\right)^{\delta_{11}}~
\mbox{ (by Lemma \ref{lem-ch-1})}.
\end{eqnarray*}

\noindent
{\rm $\mathbf{Step~ 3.} $} If $1-2\sigma<r$, then, by \cite[Theorem
2]{CP}, we conclude that there are  constants $C_{34}>0$ and $\delta_{12}\in(0,1)$ such that
$$|f(\zeta_{1})-f(\zeta_{2})|\leq
2^{\delta_{12}}C_{34}d_{\Omega}(f(z))\left(\frac{|\zeta_{1}-\zeta_{2}|}{1-|z|}\right)^{\delta_{12}}.
$$
The proof of this theorem is complete. \hfill $\Box$

\vspace{6pt}

The following result is an improvement of \cite[Lemma 3]{CP}.

\begin{lem}\label{lem-1}
%For $K\geq1$, let $f\in{\mathcal S}_{H}$ be a $K$-quasiconformal
%harmonic mapping from $\mathbb{D}$ onto a bounded domain $G$.
Let $f\in{\mathcal S}_{H}(K,\Omega)$, where $G=f(\ID)$ is a bounded domain. If
there are constants $C_{35}>0$ and $\delta_{13}\in(0,1)$ such that
for each $\zeta\in\partial\mathbb{D}$ and for $0\leq r\leq\rho<1,$
\be\label{eq-14}
\|D_{f}(\rho\zeta)\|\leq C_{35}
\|D_{f}(r\zeta)\|\left(\frac{1-\rho}{1-r}\right)^{\delta_{13}-1},
\ee
then, for $a\in\mathbb{D}$, we have
$$\diam f(I(a))\leq 32KC_{36}d_{G}(a), \quad C_{36}=2\pi e^{(1+\alpha)\pi} +\frac{2C_{35}e^{(1+\alpha)\pi}+C_{35}}{\delta_{13}},
$$
where $I(a)=\{z\in\partial\mathbb{D}:\, |\arg z-\arg a|\leq\pi(1-|a|)\}.$
%$$C_{36}=2e^{(1+\alpha)\pi} +\frac{2C_{35}e^{(1+\alpha)\pi}+C_{35}}{\delta_{13}}
%$$
%and  $\alpha$ is the same as in Theorem {\rm \ref{thm-c1}.}

\end{lem}
\bpf For $a\in\mathbb{D}$, let $a=\rho\zeta$ with $\rho=|a|$. For $z\in I(a)$, by Lemma \Ref{Lem-CP}, we have
$$%\be\label{eq-15}
|f(z\rho)-f(\rho\zeta)|\leq\int_{\gamma'}\rho\|D_{f}(\rho\xi)\| \,|d\xi|\leq 2e^{(1+\alpha)\pi}\rho \|D_{f}(\rho\zeta)\| \ell(\gamma'),
$$
where $\gamma'$ is the smaller subarc of $\partial\mathbb{D}_{\rho}$ between $\rho z$ and $\rho\zeta$, so that
$$\ell(\gamma') =\int_{\gamma'} |d\xi|=\rho |\arg (\rho\zeta)-\arg z|\leq \pi\rho (1-\rho)\leq \pi(1-\rho).
$$
Therefore,
\be\label{eq-15} |f(z\rho)-f(\rho\zeta)|\leq 2\pi e^{(1+\alpha)\pi}(1-\rho)\|D_{f}(\rho\zeta)\|.
\ee
%
%\beq\label{eq-15}
%|f(z\rho)-f(\rho\zeta)|&\leq&\int_{\gamma'}\rho\|D_{f}(\rho\xi)\| \,|d\xi|\\
%\nonumber
%&\leq&2e^{(1+\alpha)\pi}\rho\int_{\gamma'}\|D_{f}(\rho\zeta)\|
%\,|d\xi| ~\mbox{ (by Lemma \Ref{Lem-CP})}\\
%\nonumber&=&2e^{(1+\alpha)\pi}\rho\ell(\gamma')\|D_{f}(\rho\zeta)\|\\
%\nonumber
%&=&2e^{(1+\alpha)\pi}\rho^{2}\|D_{f}(\rho\zeta)\| \,|\arg (\rho\zeta)-\arg z|\\
%\nonumber&\leq&2e^{(1+\alpha)\pi}\rho^{2}(1-\rho)\|D_{f}(\rho\zeta)\|\\
%\nonumber &\leq&2e^{(1+\alpha)\pi}(1-\rho)\|D_{f}(\rho\zeta)\|,
%\eeq
Next, we have
\beq\label{eq-16}
\nonumber|f(z\rho)-f(z)|&\leq&\int_{\rho}^{1}\|D_{f}(tz)\|\,dt  \\
\nonumber
&\leq&C_{35}\int_{\rho}^{1}\|D_{f}(\rho z)\|\left(\frac{1-t}{1-\rho}\right)^{\delta_{13}-1}dt ~\mbox{ (by (\ref{eq-14}))}\\
\nonumber &=&\frac{C_{35}}{\delta_{13}}(1-\rho)\|D_{f}(\rho z)\|\\
 &\leq&\frac{2C_{35}e^{(1+\alpha)\pi}}{\delta_{13}}(1-\rho)\|D_{f}(\rho
\zeta)\|~\mbox{ (by Lemma \Ref{Lem-CP})}
\eeq
and, finally,
\beq\label{eq-17}
\nonumber |f(\zeta\rho)-f(\zeta)|&\leq&\int_{\rho}^{1}\|D_{f}(t\zeta)\|\,dt  \\
\nonumber
&\leq&C_{35}\int_{\rho}^{1}\|D_{f}(\rho\zeta)\|\left(\frac{1-t}{1-\rho}\right)^{\delta_{13}-1}\,dt ~\mbox{ (by (\ref{eq-14}))}\\
  &=&\frac{C_{35}}{\delta_{13}}(1-\rho)\|D_{f}(\rho\zeta)\|.
\eeq
%where $\gamma'$ is the smaller subarc of $\partial\mathbb{D}_{\rho}$ between $\rho z$ and $\rho\zeta$.
Again, for $z\in I(a)$, by  (\ref{eq-15}), (\ref{eq-16}), (\ref{eq-17}) and the triangle inequality, we obtain
\beq%\label{eq-20}
\nonumber |f(\zeta)-f(z)|&\leq& |f(\rho\zeta)-f(\rho z)|+|f(z)-f(\rho z)|+|f(\rho\zeta)-f(\zeta)|\\
\nonumber &\leq&C_{36}(1-\rho)\|D_{f}(\rho\zeta)\|\\
\nonumber &\leq&16KC_{36}d_{G}(a)  ~\mbox{ (by Lemma \ref{lem-ch-1})},
\eeq
which in turn implies that
$\diam f(I(a))\leq 32KC_{36}d_{G}(a)$ and the proof of the lemma is complete.
\epf

%The following result can be found in \cite{Small}.

 For $p\in(0,\infty]$, the {\it generalized Hardy space
$H^{p}_{g}(\mathbb{D})$} consists of all those functions $f:\,\mathbb{D}\rightarrow\mathbb{C}$ such that $f$ is measurable,
$M_{p}(r,f)$ exists for all $r\in(0,1)$ and  $ \|f\|_{p}<\infty$,
where
$$\|f\|_{p}=
\begin{cases}
\displaystyle\sup_{0<r<1}M_{p}(r,f)
& \mbox{if } p\in(0,\infty)\\
\displaystyle\sup_{z\in\mathbb{D}}|f(z)| &\mbox{if } p=\infty
\end{cases},
$$
and %~\mbox{ and }~
$$
M_{p}^{p}(r,f)=\frac{1}{2\pi}\int_{0}^{2\pi}|f(re^{i\theta})|^{p}\,d\theta.
$$
We refer to \cite{CPR} for more details on $H^{p}_{g}(\mathbb{D})$.

\subsection*{Proof of Theorem \ref{thm-1}}
{\rm $\mathbf{Case~ 1.} $} Let $f=h+\overline{g}\in{\mathcal S}_{H_{2}}(K, G)$, where $G$ is a bounded Pommerenke interior domain.
 Then, by definition, \eqref{CP-extraeq6} holds and thus (see for example, \cite[Proof of Theorem 3]{Po}),
%$$\sup_{z\in\mathbb{D}}\left|(1-|z|^{2})\frac{h''(z)}{h'(z)}-2\overline{z}\right|<4,
%$$
%which implies that
there are constants $\rho_{0}\in(0,1)$ and $\beta_{1}>0$ such that, for $\rho_{0}\leq\rho<1$ and
$\theta\in[0,2\pi]$,
\be\label{eq-ch-5}
\mbox{Re}\left[e^{i\theta}\frac{h''(\rho e^{i\theta})}{h'(\rho e^{i\theta})}\right]
\geq-\frac{1-\beta_{1}}{1-\rho}.
\ee
For $\rho_{0}\leq r\leq\rho<1$,  by integrating both sides of (\ref{eq-ch-5}), we have
$$(1-r)^{\beta_{1}-1}|h'(re^{i\theta})|\leq(1-\rho)^{\beta_{1}-1}|h'(\rho e^{i\theta})|,
$$
which, by \eqref{CP-extraeq7}, deduces that
\be\label{eq-ch-6}
(1-r)^{\beta_{1}-1}\|D_{f}(re^{i\theta})\|\leq\frac{2K}{1+K}(1-\rho)^{\beta_{1}-1}\|D_{f}(\rho e^{i\theta})\|.
\ee
%{\rm $\mathbf{Step~ 2.} $}
For $\rho\in(\rho_{0},1)$, we choose a positive integer $N$ and $r_{0},\ldots,r_{N}$ with
$r_{N}=\rho_{0}<r_{N-1}<\cdots<r_{1}<r_{0}=\rho$ such that, for $n\in\{0,1,\ldots, N-1\}$,
$$2^{n}(1-\rho)\leq1-r_{n}<2^{n+1}(1-\rho).
$$
For $\theta\in[0,2\pi),$ let
$$I(r_{n}e^{i\theta})=\{\zeta\in\partial\mathbb{D}:\,|\arg\zeta-\theta|\leq\pi(1-r_{n})\}.
$$
For $2\leq n\leq N$ and $e^{it}\in I(r_{n}e^{i\theta})\backslash
I(r_{n-1}e^{i\theta})$, let $\varphi=t-\theta$. Then, for $2\leq
n\leq N$,
\be\label{eq-che}
\pi(1-r_{n-1})\leq|\varphi|\leq\pi(1-r_{n}).
\ee

By the assumption, we let
\be\label{eq-p-5}
c_{p}=\sup_{0<r<1}\left\{\sup_{w_{1},w_{2}\in\gamma_{r}}\frac{\ell\big(\gamma_{r}
[w_{1},w_{2}]\big)}{d_{G_{r}}(w_{1},w_{2})}\right\}<\infty,
\ee
where $\gamma_{r}$ is given by (\ref{eq-y}). Then,  by (\ref{eqy}),
(\ref{eq-p-5}) and \cite[Theorem 4]{CP},  $\|D_{f}\|\in
H^{1}_{g}(\mathbb{D}).$ Hence, for  $n\in\{0,1,\ldots, N-1\}$, by
(\ref{eqy}), (\ref{eq-p-5}), Lemma \ref{lem-1} and \cite[Inequality
(2.3)]{CP}, there is a positive constant $C_{37}$ such that

%\textcolor[rgb]{1.00,0.00,0.00}{
\beq\label{eq-p-4}
\nonumber\frac{1}{K}\int_{I(r_{n}e^{i\theta})}\|D_{f}(
e^{it})\|dt&\leq&\int_{I(r_{n}e^{i\theta})}l(D_{f}( e^{it}))\,dt\\
\nonumber &\leq&\int_{I(r_{n}e^{i\theta})}|df( e^{it})| =\ell(I(r_{n}e^{i\theta}))\\
%\nonumber &=&\ell(I(r_{n}e^{i\theta})) \\
\nonumber&\leq&c_{p}\diam\big(f(I(r_{n}e^{i\theta}))\big)~\mbox{ (by (\ref{eq-p-5}))}\\
\nonumber &\leq&C_{37}c_{p}d_{G}(r_{n}e^{i\theta})~\mbox{ (by (\ref{eqy}) and Lemma \ref{lem-1})}\\
&\leq&\frac{2KC_{37}c_{p}}{1+K}(1-r_{n})\|D_{f}(r_{n}e^{i\theta})\|
~\mbox{ (by \cite[Inequality (2.3)]{CP}).}
\eeq
Let $I_{n}(\theta)=I(r_{n}e^{i\theta}).$ Since
$\partial\mathbb{D}=I_{0}(\theta)\cup(I_{1}(\theta)\setminus
I_{0}(\theta))\cap\cdots \cap (I_{N}(\theta)\setminus
I_{N-1}(\theta)),$ by (\ref{eq-ch-6}) and (\ref{eq-che}), we see
that
\beq\label{CP-extraeq8}
\Lambda_{f}&=&\int_{0}^{2\pi}\|D_{f}(e^{it})\|\frac{1-\rho^{2}}{|e^{it}-\rho e^{i\theta}|^{2}}\,dt \leq J_0 +\sum_{n=1}^{N}J_n
%=\int_{0}^{2\pi}\|D_{f}(e^{it})\|\frac{1-\rho^{2}}{|1-\rho e^{i(\theta-t)}|^{2}}\,dt
%\\
%\nonumber&\leq& J_0 +\sum_{n=1}^{N}J_n
\eeq
where, by (\ref{eq-p-4}),
\be\label{CP-extraeq9}
J_0=\frac{2}{1-\rho}\int_{I_{0}(\theta)}\|D_{f}(e^{it})\|\,dt \leq
\frac{4K^{2}C_{37}c_{p}}{1+K}\|D_{f}(\rho e^{it})\|
%, ~\mbox{(by(\ref{eq-p-4}))}
\ee
and
\beq\label{eq-p-8} \nonumber J_n&=&
\int_{I_{n}(\theta)\setminus
I_{n-1}(\theta)}\|D_{f}(e^{it})\|\frac{1-\rho^{2}}{|1-\rho e^{i(\theta-t)}|^{2}}\,dt\\
\nonumber
&=& \int_{I_{n}(\theta)\setminus
I_{n-1}(\theta)}\frac{\|D_{f}(e^{it})\|(1-\rho^{2})}{(1-\rho)^{2}+4\rho\sin^{2}\frac{\varphi}{2}}\,dt \\
\nonumber
&\leq& \int_{I_{n}(\theta)\setminus
I_{n-1}(\theta)}\|D_{f}(e^{it})\|\frac{\pi^{2}(1-\rho^{2})}{4\rho\varphi^{2}}\,dt\\
 \nonumber
&\leq& \int_{I_{n}(\theta)\setminus I_{n-1}(\theta)}\frac{\|D_{f}(e^{it})\|(1-\rho^{2})}{4\rho(1-r_{n-1})^{2}}\,dt~\mbox{ (by (\ref{eq-che}))}\\
\nonumber
&\leq& \int_{I_{n}(\theta)}\frac{\|D_{f}(e^{it})\|(1-\rho^{2})}{4\rho(1-r_{n-1})^{2}}\,dt\\
\nonumber
&\leq& \frac{K^{2}C_{37}c_{p}}{1+K} \left (\frac{(1-\rho)(1-r_{n})\|D_{f}(r_{n}e^{it})\|}{\rho(1-r_{n-1})^{2}}\right )~\mbox{(by
(\ref{eq-p-4}))}
\\ \nonumber
&\leq& \frac{16K^{2}C_{37}c_{p}}{(1+K)\rho}\left (
\frac{(1-\rho)\|D_{f}(r_{n}e^{it})\|}{1-r_{n}}\right ). \eeq
By (\ref{eq-ch-6}), \eqref{CP-extraeq8}, \eqref{CP-extraeq9} and the last inequality, we conclude that
\beq
\nonumber \Lambda_{f} &\leq
&\frac{4K^{2}C_{37}c_{p}}{1+K}\|D_{f}(\rho
e^{it})\|\left(1+\frac{4}{\rho}\sum_{n=1}^{N}\frac{(1-\rho)}{(1-r_{n})}\frac{\|D_{f}(r_{n}e^{it})\|}{\|D_{f}(\rho
e^{it})\|}\right)\\ \nonumber %\eeq \beq\nonumber
&\leq&\frac{4K^{2}C_{37}c_{p}}{1+K}\|D_{f}(\rho
e^{it})\|\left(1+\frac{8K}{1+K}\frac{1}{\rho}\sum_{n=1}^{N}\frac{1}{2^{n\beta_{1}}}\right)~\mbox{ (by
(\ref{eq-ch-6}))}.
\eeq
Thus,
$$\sup_{\zeta\in\mathbb{D}\backslash\overline{\mathbb{D}}_{\rho_{0}}}\frac{1}{2\pi}
\int_{\partial\mathbb{D}}\frac{\|D_{f}(\xi)\|}{\|D_{f}(\zeta)\|}
\frac{1-|\zeta|^{2}}{|\xi-\zeta|^{2}}|d\xi|<\infty,
$$
%where
%$$\Lambda_{f}=\int_{0}^{2\pi}\|D_{f}(e^{it})\|\frac{1-\rho^{2}}{|e^{it}-\rho e^{i\theta}|^{2}}\,dt
%$$
%and $I_{n}(\theta)=I(r_{n}e^{i\theta})$.

{\rm $\mathbf{Case~ 2.} $} For $\rho\in[0,\rho_{0}]$ and
$\theta\in[0,2\pi]$, by \cite[Theorem 4]{CP}, we have
%\be\label{eq-ch-8}
$$\int_{\mathbb{D}}\|D_{f}(\xi)\|\frac{(1-\rho^{2})}{|\xi-\rho
e^{i\theta}|^{2}}|d\xi|\leq\frac{1+\rho_{0}}{1-\rho_{0}}\int_{\mathbb{D}}\|D_{f}(\xi)\|\,|d\xi|<\infty.
$$
On the other hand, for $\theta\in[0,2\pi]$,
%\be\label{eq-ch-9}
$$\min_{\rho\in[0,\rho_{0}]}\|D_{f}(\rho e^{i\theta})\|>0.
$$
In this case, Theorem \ref{thm-1} follows from the last two inequalities.
% (\ref{eq-ch-8}) and (\ref{eq-ch-9}).
%The proof of the theorem is complete.
\hfill $\Box$

%\hfill $\Box$
%\begin{Lem}\label{Lem-D}
%Let $f=h+\overline{g}\in{\mathcal S}_{H}$, where $h$ and $g$ are
%analytic in $\mathbb{D}$. Then, for $z\in\mathbb{D}$,
%$$\left|(1-|z|^{2})\frac{h''(z)}{h'(z)}-2\overline{z}\right|\leq2\alpha,$$
%where  $\alpha$ is the same as in the Lemma {\rm\Ref{Lem-CP}}.

%\end{Lem}

%\subsection*{Proof of Theorem \ref{thm-1.02}}

%\hfill $\Box$

%\subsection*{Proof of Theorem \ref{thm-1.0}}
 %\hfill $\Box$

%\bigskip

%{\bf Acknowledgements:} We are grateful to the referee for his/her
%comments and suggestions. This research was partly supported by the
%National Natural Science Foundation of China (No. 11401184 and No.
%11571216), the Hunan Province Natural Science Foundation of China
%(No. 2015JJ3025), the Excellent Doctoral Dissertation of Special
%Foundation of Hunan Province (higher education 2050205), the
%Construct Program of the Key Discipline in Hunan Province. The
%second author is currently on leave from
 % Indian Institute of Technology Madras.

%{\bf Acknowledgements:} This research was  partly  supported by NSF
%of China (No. 11071063). This work was also supported in part by the
%Construct Program of the Key Discipline in Hunan Province (No.
%[2011] 76) and the Start Project of Hengyang Normal University (No.
%12B34).
%Acta Math. Sinica (English Series) 32 (3), 297-308

\normalsize

\end{document}